\documentclass[preprint,number]{elsarticle}
\usepackage{bm}
\usepackage{amssymb,amsmath,amsfonts}
\usepackage{mathrsfs}
\usepackage{mathabx}
\usepackage{subfigure}
\usepackage{afterpage}
\usepackage{color}
\usepackage{tabu}

\renewcommand{\vec}[1]{\mbox{\boldmath$#1$}}
\newcommand{\dif}{\mathrm{d}}

 \title{Benchmark computations for the polarization tensor characterization of  small conducting objects}

 \author[1]{A.A.S. Amad}
\author[2]{P.D. Ledger\corref{cor1}}
\ead{p.d.ledger@keele.ac.uk}
\author[3]{T. Betcke}
\author[4]{D. Praetorius}
\cortext[cor1]{Corresponding author}
\affiliation[1]{organization={Zienkiewicz Centre for Computational Engineering, College of Engineering, Swansea University Bay Campus}, addressline={Fabian Way},
postcode={SA1 8EN}, city={Swansea}, country={UK}}
\affiliation[2]{organization={School of Computing \& Mathematics, Keele University}, city={Newcastle under Lyme, Staffordshire},
           postcode={ST5 5BG},
           country={UK}}
  \affiliation[3]{organization={Centre for Inverse Problems, University College London}, addressline={Gower Street}, city={London},
           postcode={WC1E 6BT},
           country={UK}} 
  \affiliation[4]{organization={TU Wien, Institute for Analysis and Scientific Computing}, addressline={Wiedner Hauptstr. 8-10 / E101 / 4}, city={Wien},
           postcode={1040},
           country={Austria}} 
                   
\date{Submitted 2nd November 2021, Revised 29th March 2022, 10th June 2022}
 
\begin{document}
 \begin{abstract}
The characterisation of small low conducting inclusions in an otherwise uniform background from low-frequency electrical field measurements has important applications in medical imaging using electrical impedance tomography as well as in geological imaging using electrical resistivity tomography. It is known that such objects can be characterised by a P\'oyla-Szeg\"o (polarizability) tensor. Such characterisations have attracted interest as they can provide object features in a machine learning classification algorithm and provide an alternative imaging solution. However, to be able train machine learning algorithms, large dictionaries are required and it is essential that the characterisations are accurate. In this work, we obtain accurate numerical approximations to the tensor coefficients, by applying an adaptive boundary element method. The goal being to provide a sequence of benchmark {computations} for the tensor coefficients to allow other software developers check the accuracy of their codes.
 \end{abstract}           
           

 \maketitle

\section{Introduction}

The purpose of this paper is provide a series of accurate benchmark {computations} for the P\'oyla-Szeg\"o  tensor (PST) characterisation of small conducting objects with low contrasts. Our benchmark {computations} are intended to be useful for developers of computational tools for solving partial differential equations and, in particular, for those developing tools that will form part of a machine learning (ML) classification algorithm for detecting and classifying small conducting inclusions.

{The characterisation of conducting objects with low conductivities in an otherwise uniform background using low frequency electric fields has a wide range of important applications including medical imaging and geology using electrical impedance tomography (EIT) and electrical resistivity tomography, respectively  e.g. \cite{LionheartPM2004}. In particular, being able to characterise small objects is important for the related inverse problem where the task is to be able to determine the shape and material contrasts of hidden conducting inclusions from measurements of the induced voltage at the boundary from induced current distributions \cite{MuellerBookSIAM2012}. This inverse problem was first studied by Calder\'on  \cite{CalderonCAM2006}. There is a large and increasing literature on the EIT inverse problem and many different computational approaches have been proposed.  Some examples of alternative solution approaches include total variational regularisation  \cite{ChungJCP2005}, topological derivative \cite{HintermullerACM2012} and machine learning approaches \cite{RymarczykSensors2018}. }

{For EIT it has been rigorously proven by Cedio-Fengya, Moskow and Vogelius~\cite{CedioIP1998} that} an asymptotic expansions of the perturbed field for small objects {simply connected}
 in the form $B_\alpha := \alpha B + {\bm z}$, where $\alpha$ denotes the object size, $B$ a unit sized object placed at the origin and ${\bm z}$ is the translation, {in a bounded domain}
 has a leading order term, 
\begin{align} \label{eq:asymptoticExpansion}
    (u_\alpha - u_0)(\bm{x})  =&  \nabla_{\bm{x}} G(\bm{x}, \bm{z})\cdot [{\mathcal T}(k, \alpha B) \nabla_{\bm{z}} u_0(\bm{z})]  + O(\alpha^{4})
\end{align} 
 as $\alpha \to 0$ where $u_\alpha$ denotes the field in the presence of the object, $u_0$ denotes the background field,
 $G(\bm{x}, \bm{z}) := 1/{(4\pi |\bm{x} - \bm{z}|)}$ is the free space Laplace Green's function, $k$ denotes the (conductivity) contrast, ${\mathcal T}(k, \alpha B)=\sum_{i,j=1}^3{\mathcal T}_{ij} {\bm e}_i \otimes {\bm e}_j$ is the appropriate polarizability tensor (also sometimes known as a polarization tensor) for this problem and ${\bm e}_i$ the orthonormal coordinate directions.  {The appropriate polarization tensor for this problem coincides with the PST, which is independent of the object's position.  {Ammari and Kang~\cite[Section 4.1]{ammarikangbook}} have derived the corresponding {complete} asymptotic expansion of this problem and have introduced the concept of generalised polarization tensors, which offer improved object characterisations and have the PST as the simplest type. The same PST object description is also known to characterise scatterers in electromagnetic scattering~\cite{ledgerlionheart2012}  and it has been proposed that electro-sensing fish use them for sensing their pray \cite{taufiqBB2016}. A closely related polarization tensor also describes conducting objects in metal detection~\cite{ledgerlionheart2014}.}

In the majority of  EIT approaches, the hidden object's position and its shape is inferred by finding information about the conductivity distribution in the imaging region. However, an alternative approach is offered by combining dictionaries of object characterisations using PSTs with classification techniques, such as those in ML  algorithms (see \cite{bishopbook2006} for a discussion of different approaches), due to the separation of object characterisation and in position in (\ref{eq:asymptoticExpansion}). {A ML approach to object classification using polarizability tensors has already been shown to be effective for metal detection~\cite{wilson2021identification}, but for such an approach} to be effective for EIT problems, a description of accurate PST coefficients is required.

In this paper, we provide novel computational {benchmark computational} solutions to a series of PST object characterisations. {The non-dimensional benchmark geometries considered are a unit radius sphere, an ellipsoid with principal semi axes $a=1$, $b=0.7$ and $c=0.5$, an L-shape domain
$ B_1 \bigcup B_2$ with  $B_1 = [0.0,7.8] \times [0.0,2.0] \times [0.0,1.5]$, $B_2 =  [0.0,2.2] \times [2.0,5.6] \times [0.0,1.5]$, the cube $[-1,1]^3$, the tetrahedron with vertices
\begin{equation}
{\vec v}_1 = \begin{pmatrix} 
0 \\
0 \\
0 
\end{pmatrix}, \ 
{\vec v}_2= \begin{pmatrix} 
7 \\
0 \\
0 
\end{pmatrix}, \  
{\vec v}_3 =  \begin{pmatrix} 
5.5 \\
4.6 \\
0.0 
\end{pmatrix}, \
{\vec v}_4=  \begin{pmatrix} 
3.3 \\
2.0 \\
5.0 
\end{pmatrix},
\nonumber
\end{equation}
and the key contained in the rectangular region $B_3  \bigcup B_4 \bigcup B_5$
where 
\begin{align*}
B_3 & = [-7.0,7.0] \times [3.0,15.0] \times [-1.25,1.25], \\
B_4 & = [-4.0,4.0] \times {[0.0,3.0]} \times [-1.25,1.25], \\
B_5 & = [-3.25, 3.25] \times [-19.0,0.0] \times [-1.25,1.25]. 
\end{align*}
Each object geometry is centered at the origin and they are scaled by different length scales $\alpha$ in meters. The benchmark computations for the PST coefficients are obtained }by post processing the approximate solutions to a transmission problem, which is governed by Laplace's equation set on an unbounded region surrounding the (scaled) object of interest. Our benchmark {computations} are not only of interest to EIT practitioners, but are also of interest to other computational partial differential equation solver developers since the PST coefficients can be obtained by a simple post-processing of transmission problem solutions and, by comparing solutions, can provide a measure of accuracy of a given scheme. To ensure our solutions are accurate, and 
inspired by a previous boundary integral computation of PST coefficients \cite{taufiqDMaths2013}, the transmission problem will be solved using a boundary integral approach using the boundary element python package (Bempp, formally known as the BEM++ library)~\cite{SmigajACM2015}, leading to approximate tensor coefficients for different objects. An adaptive boundary integral scheme, which has already been shown to be effective for a capacitive problem~\cite{BetckeJCP2019}, is novelly applied to the computation of the PST to generate a set of accurate benchmark {computations} for different geometries.

The work proceeds as follows: In Section \ref{sec:ModelProblem}, the mathematical model and boundary integral/layer potential formulae to compute the PST are introduced. 
Then, in Section \ref{sec:BEM:discretisation}, a brief description of how  Bempp  can be applied to compute the PST is presented and, in Section \ref{sec:AdaptiveMesh}, the application of an existing adaptive algorithm to the computation of the PST is described. Section \ref{sec:numericalresults} describes {the geometries used for our benchmark computations} and presents numerical results for the accurate computation of the PST using the adaptive algorithm. The paper ends with some concluding remarks in Section \ref{sec:Conclusions}.

\section{Small Object Characterisation} \label{sec:ModelProblem}

\subsection{Model Problem}
\begin{figure}[h]
\centering
\includegraphics[scale=0.5]{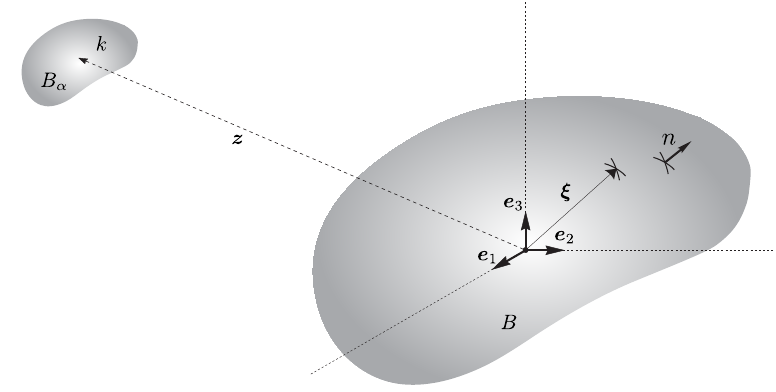} 
  \caption{Illustration of physical object $B_\alpha$ described as $B_\alpha =\alpha B + {\bm z}$. }
        \label{fig:physobject}
\end{figure}

{We following the formulation in Cedio-Fengya, Moskow and Vogelius~\cite{CedioIP1998} and Ammari and Kang~\cite[Section 4.1]{ammarikangbook} restricted to the characterisation of a single  Lipschitz bounded simply connected region $B_\alpha$ in $\Omega \subset {\mathbb R}^3$. We let $0<k\ne1< + \infty $ be the conductivity (contrast) of $B$ to the background (i.e $k=\gamma |_{B_\alpha}/ \gamma |_{\Omega \setminus B_\alpha}$ where $\gamma $ is the electrical conductivity). We  set $B_\alpha := \alpha B + {\bm z}$, where $\alpha\ll 1$ denotes the (small) object size parameter, $B$ a unit sized object placed at the origin, and ${\bm z}$ is the translation, as illustrated in Figure~\ref{fig:physobject}. The voltage $U$ in absence of the inclusion is given by the solution to
\begin{align*}
\nabla^2  U & = 0 && \text{in $\Omega $}\\
\frac{\partial U}{\partial {\vec n}} & = \psi && \text{on $\partial \Omega$} 
\end{align*}
where $\psi$ denotes the imposed current such that $\int_{\partial \Omega} \psi \dif {\vec x} =0$ and ${\vec n}$ is the unit outward normal to $\partial \Omega$. In the presence of the inclusion, the voltage $u_\alpha$ is given by the solution to
\begin{align*}
\nabla \cdot (( 1 + (k-1)\chi(B) ) \nabla u_\alpha )  & = 0 && \text{in $\Omega $}\\
\frac{\partial u_\alpha}{\partial {\vec n}} & = \psi && \text{on $\partial \Omega$}
\end{align*}
For this problem, Cedio-Fengya {\em et al} prove an asymptotic expansion of the form stated in (\ref{eq:asymptoticExpansion}) where the real symmetric  rank 2 PST ${\mathcal T} := {\mathcal T}(k, \alpha B)$ is a function of  $k$,  $\alpha $, and $B$, but it is independent of ${\vec z}$. Ammari and Kang obtained similar asymptotic expansions for other closely related problems in EIT. Both Cedio-Fengya {\em et al} and  Ammari and Kang discuss the extension to multiple inclusions in EIT. Properties of the tensor ${\mathcal T}$ are discussed by Cedio-Fengya {\em et al}  and Ammari and Kang, which we will draw on below. Ammari and Kang also discuss algorithms for recovering PST coefficients from field measurements.}

{Our interest in this work focus on computational schemes for accurately computing the coefficients of the real symmetric  rank 2 PST ${\mathcal T} $ for known $\alpha$, $B$, and $k$. We start by recalling a series of alternative, but equivalent formulations of  ${\mathcal T}$.}

\subsection{Transmission Problem Formulation}
Given the solutions $\vartheta_i$, $i=1,2,3$, to the transmission problem
\begin{subequations}\label{eq:TransmissionProblem}
    \begin{align}  \label{eq:TransmissionProblem:a}
        \nabla^2 \vartheta_i  & = 0 & & \hbox{in $B \bigcup B^c$}  ,  \\
        \left [ \vartheta_i \right ]_\Gamma  & = 0 && \hbox{on $\Gamma$},  \label{eq:TransmissionProblem:b}\\
        {\vec n} \cdot \nabla \vartheta_i \Big|_+  - {\vec n} \cdot k \nabla \vartheta_i \Big|_- & = {\vec n} \cdot \nabla_{\bm{\xi}} \xi_i && \hbox{on $\Gamma$},  \label{eq:TransmissionProblem:c}\\
        \vartheta_i   &\to 0 && \hbox{as $|\bm{\xi}| \to \infty$} , 
    \end{align} 
\end{subequations}
where $B^c := \mathbb{R}^3\setminus \overline{B} $, $\Gamma:= \partial B$,  $\left [ \cdot \right ]_\Gamma= \cdot |_+ - \cdot|_-$ denotes the jump,  $\big|_-$ the interior trace and $\big|_+$ the exterior trace,  $\xi_j = {\bm e}_j \cdot {\bm \xi}$ denotes the $j$-th the component of ${\bm \xi} $ measured from an origin in $B$ and differentiation is with respect to ${\bm \xi}$, the coefficients of the PST %
\begin{align}\label{eq:PStensor}
    {\mathcal T}_{ij} :=  \alpha^3 \left ( \left (k - 1 \right) |B|\delta_{ij} + \left (k - 1 \right)^2 \int_\Gamma {\vec n} \cdot  \nabla \vartheta_i  \Big|_- \xi_j\, \dif \bm{\xi} \right ),
\end{align}
follow by {simplifying the result for the coefficients of a GPT in~\cite[equation (4.10)]{ammarikangbook}.}
In the above, $|B|$ is the volume of the object and $\delta_{ij}$ is the Kronecker delta. By integration by parts,  (\ref{eq:PStensor}) can be shown to be equivalent to
\begin{align*}
    {\mathcal T}_{ij} := & \alpha^3 \left ( \left (k - 1 \right) |B|\delta_{ij} - \left (k - 1 \right)^2 \left [  \int_{B^c} \nabla \vartheta_i \cdot \nabla \vartheta_j \dif {\bm \xi} \right . \right . \nonumber \\
    & \left . \left . + \int_B k \nabla \vartheta_i \cdot \nabla \vartheta_j \dif {\bm \xi}  \right ] \right ), 
\end{align*}
which makes the  symmetry  of the PST clear. Its independence of ${\bm z}$ can similarly be established~\cite{CedioIP1998}.

\subsection{Boundary Integral (BI) Formulation} \label{sec:BI}

By following \cite[equations (56)-(58)]{CedioIP1998} it is possible to rewrite \eqref{eq:PStensor}, as follows
\begin{align}\label{eq:BIformula}
    {\mathcal T}_{ij} =  \alpha^3 \left ( \dfrac{1}{r} |B|\delta_{ij} +  \left (1 - \dfrac{1}{r} \right) \int_\Gamma \psi_i  n_j \, \dif \bm{\xi} \right ),
\end{align}
where $r = \dfrac{1}{k}$ and $\psi_i := -k\vartheta_i$ can be expressed {as the solution to the boundary integral equation}~\cite{CedioIP1998}
\begin{equation}\label{eq:BEM:BI}
\psi_i(\bm{\xi}) := \frac{1}{r-1}\mathcal{S} \left[ \left( \lambda I + \mathcal{K}^{\ast} \right)^{-1}\left ( {\vec n}_{\bm{\nu}} \cdot \nabla_{\bm{\nu}} \nu_i \right ) \right] (\bm{\xi}), \quad  \forall \bm{\xi} \in \mathbb{R}^3.
\end{equation}
Here, $\mathcal{S}$ is the single layer potential operator~{\cite[equation (2.12)]{ammarikangbook},} 
\begin{equation}
        \mathcal{S}[f]({\vec \xi}) := \frac{1}{4 \pi } \int_{\Gamma} \frac{1}{|\bm{\xi} - \bm{\nu}|} f({\vec \nu})\, \dif {\vec \nu}, \nonumber
        \end{equation}
{$4\pi $ is the area of the unit sphere in $\mathbb{R}^3$}, $I$ is the identity operator and $\mathcal{K}^*$ is the adjoint double layer boundary operator \cite[equation (2.20)]{ammarikangbook},
\begin{equation}
    \mathcal{K}^{\ast} [f](\bm{\xi}) := \frac{1}{4\pi } \mbox{ p.v.}\int_{\Gamma} \frac{ (\bm{\xi} - \bm{\nu})\cdot {\vec n}_{\bm{\xi}} }{|\bm{\xi} - \bm{\nu}|^3}  f(\bm{\nu})\, \dif \bm{\nu}. \nonumber
\end{equation}
In addition,  p.v. denotes the Cauchy principal value, and $|\bm{\xi} - \bm{\nu}|$ is the distance between $\bm{\xi}$ and $\bm{\nu}$, for further details see~{\cite[page 18]{ammarikangbook}.}

\subsection{Layer Potential (LP) Formulation} \label{sec:LayerPotential}

In~{\cite[Section 4.1]{ammarikangbook}}, the coefficients of ${\mathcal T}$ are shown to be equivalently obtained as follows
\begin{equation}\label{eq:PSintegral}
  {\mathcal T}_{ij} = \alpha^3 \int_{\Gamma} \xi_j \phi_i\,\dif \bm{\xi},
\end{equation}
where $\phi_i$ is {the solution to the boundary integral equation}
\begin{equation}\label{eq:BEM:LP}
 \phi_i(\bm{\xi}) :=  \left (\lambda I + \mathcal{K}^* \right )^{-1} \left ( {\vec n}_{\bm{\nu}} \cdot \nabla_{\bm{\nu}} \nu_i \right )(\bm{\xi}),\quad \bm{\xi} \in \Gamma, 
\end{equation}
and $\lambda := \dfrac{k+1}{2(k-1)}$. Note that, for exact computations, we have $\psi_i(\bm{\xi}) = \dfrac{1}{r-1}\mathcal{S} \left [ \phi_i \right ] ({\bm \xi})$.

\subsection{Weighted (Alternative) Formulation} \label{sec:Alternative}

By taking a weighted average of \eqref{eq:BIformula} and \eqref{eq:PSintegral} we get
\begin{align}\label{eq:AlternativeFormula}
    {\mathcal T}_{ij} =  \alpha^3 \left (  \dfrac{\beta}{r}|B|\delta_{ij} + \int_\Gamma  \left( \beta  \left (1 - \dfrac{1}{r} \right) \psi_j  n_i + \left (1 - \beta \right)\phi_j\xi_i \right) \, \dif {\vec \xi} \right ),
\end{align}
with $0 \leq \beta \leq 1$.
Note that, for $\beta = 0$ we obtain the formula \eqref{eq:PSintegral}, whereas for $\beta = 1$ we have \eqref{eq:BIformula}.

\section{Computation of the P\'olya-Szeg\"o Tensor Coefficients} 
\label{sec:BEM:discretisation}

In this work, we use Bempp \cite{SmigajACM2015} to compute the coefficients of the PST. Bempp is a C++ library with Python bindings for the solution of boundary integral equations using Galerkin discretisations described in{~\cite{SmigajACM2015} and further improved in~\cite{BetckeACM2020}.
Note that corners,} {in contrast to collocation,} {are not an issue for Galerkin boundary element methods on piecewise smooth Lipschitz domains. We are using fully numerical Erichsen/Sauter
type singular quadrature rules that are well suited for this type of discretisation~\cite{Erichsen1998EfficientAQ}.} Given the similarity of the transmission problem
(\ref{eq:TransmissionProblem}) and the boundary integral solution formulations (\ref{eq:BEM:BI}) and (\ref{eq:BEM:LP})  to the capacitive problem already successfully considered in~\cite{BetckeJCP2019}, we will focus on a dual mesh discretisation although a primal discretisation is also possible and may also offer further efficiencies.
In practice, given a surface mesh discretisation of $B$ of size $h$ we solve for $\phi_i^h$ and obtain $\psi_i^h({\bm \xi}) = \dfrac{1}{r-1}\mathcal{S} \left [ \phi_i^h \right ] ({\bm \xi})$.
 Following the numerical solution, the discrete approximation ${\mathcal T}_{ij}^h$ to ${\mathcal T}_{ij}$ then follows from a simple postprocessing  of the solutions using (\ref{eq:AlternativeFormula}).

Although the coefficients of the true tensor ${\mathcal T}$ are symmetric, the coefficients of its numerical approximation ${\mathcal T}^h$ using boundary elements are not guaranteed to have this property, {but tend to a symmetric tensor as the mesh is refined}. Thus, to ensure the symmetry of  the computed PST,  we will symmetrise it as follows:
    \begin{equation}
         {\mathcal T}^{h,sym}_{ij} := \dfrac{1}{2} \left({\mathcal T}^h_{ij} +  {\mathcal T}^h_{ji}  \right) . \nonumber
    \end{equation}    
Henceforth, when presenting results, we will only consider the symmetrised tensor and drop the superscript $sym$. In the next section, we illustrate the convergence behaviour of the numerical scheme for uniform grid refinement for a PST with smooth analytical solution.

\subsection{Numerical Convergence Study for Ellipsoidal Objects} \label{sec:example:sphere}

\begin{figure}[h]
\centering
    \subfigure[Sphere]{\includegraphics[scale=0.4]{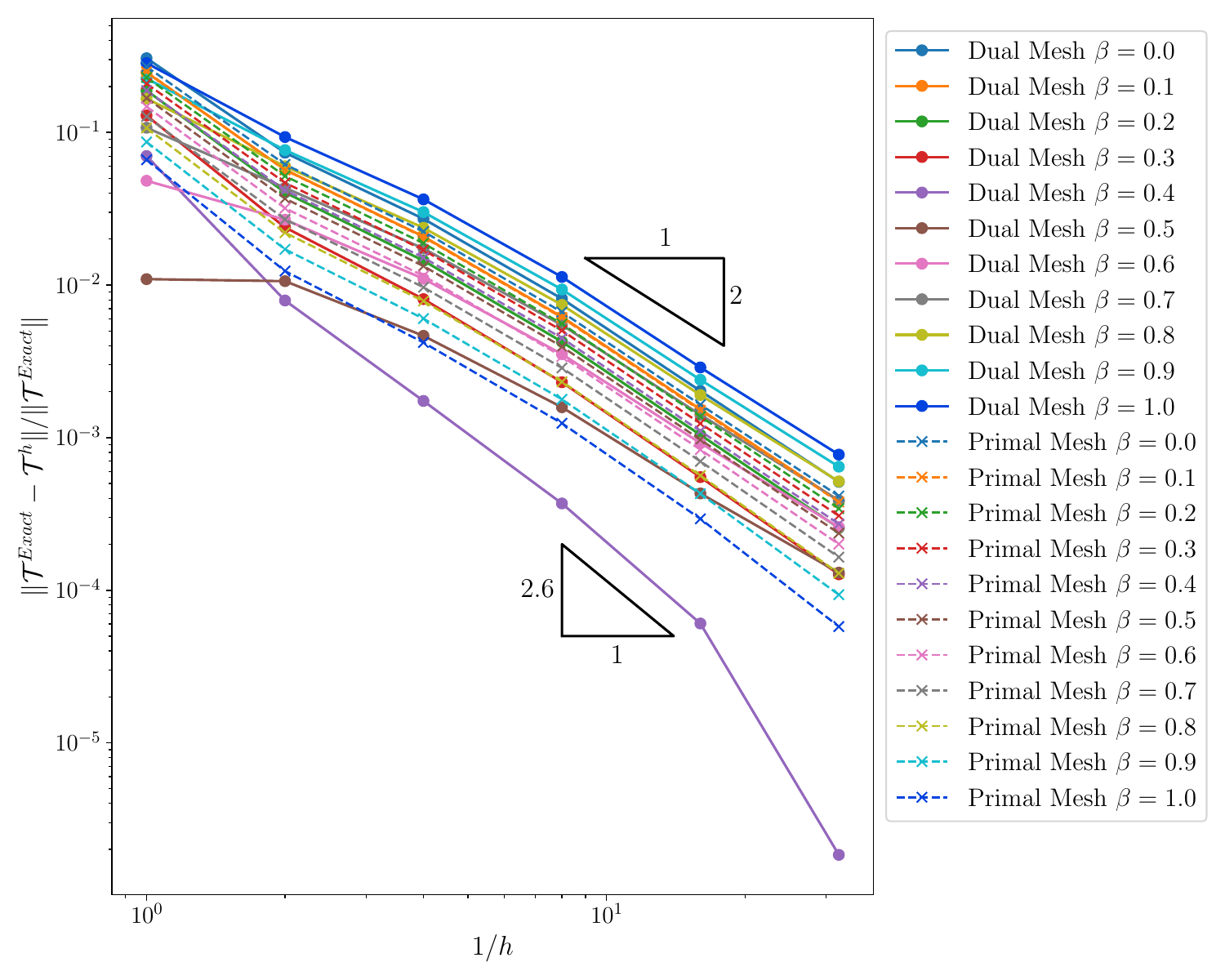} 
    \label{fig:sphere:relativeError}}
    \subfigure[Ellipsoid]{\includegraphics[scale=0.4]{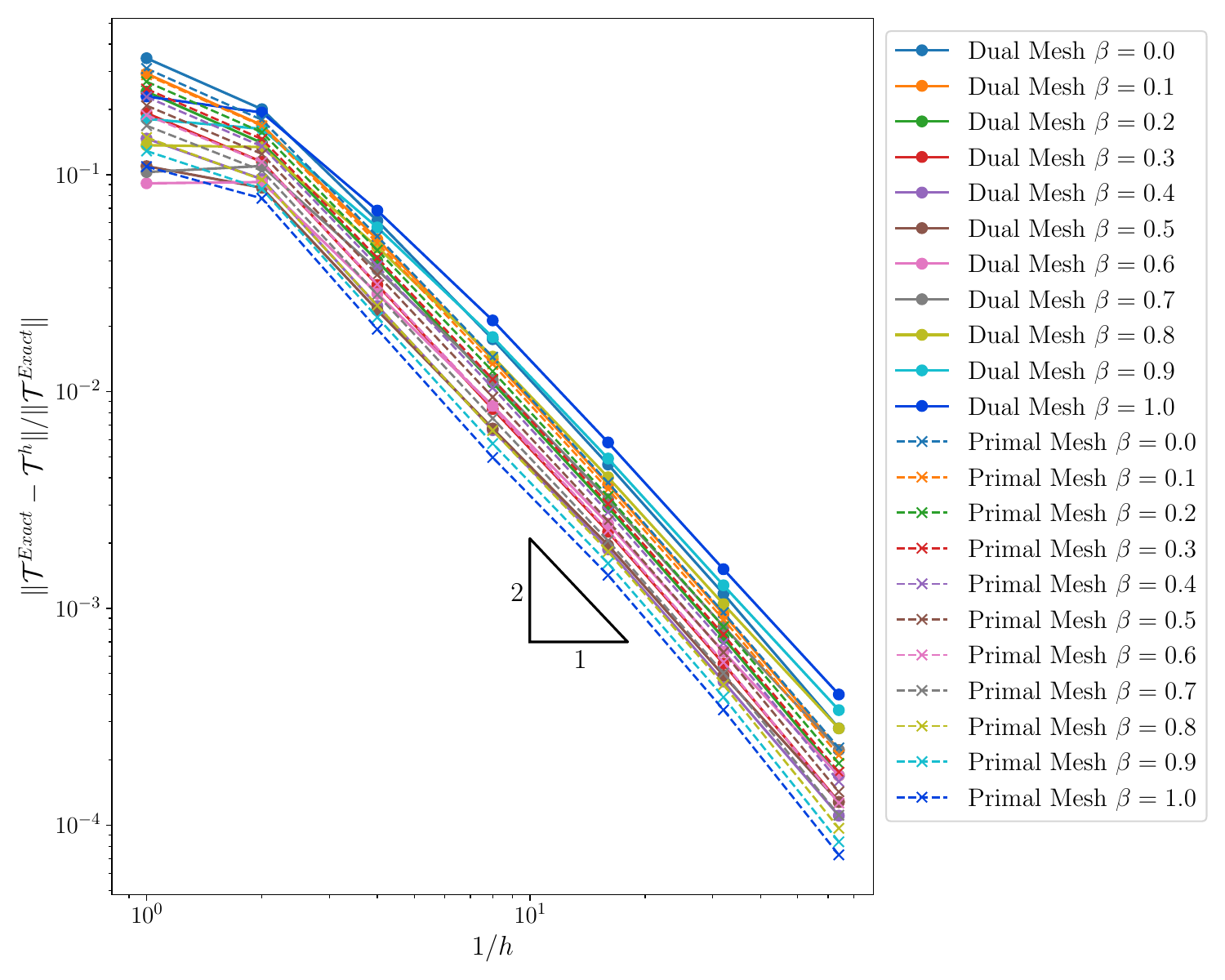}
     \label{fig:ellipsoid:relativeError}}
  \caption{Relative error between the study of the polarization tensor ${\mathcal T}^h$ \eqref{eq:AlternativeFormula} computed on primal mesh $\mathcal{M}$ and dual mesh $\mathcal{M}^{dual}$  for $\beta = 0.0, \ldots, 1.0$ and compared with the exact solution in \cite[page 83]{ammarikangbook} (a) sphere with radius $a = b = c = 1$ and (b) ellipsoid with radii $a = 1$, $ b = 0.7 $, $c = 0.5$. }
        \label{fig:benchmark:Error}
\end{figure}

The convergence of the approximate tensor $\mathcal{T}^{h}$  to the true PST $\mathcal{T}$ obtained  
 on a sequence of uniformly meshes for ellipsoidal objects is investigated. In these experiments, we choose $B_\alpha = \alpha B$ and define $B$ by
\begin{equation}
\dfrac{\xi_1^2}{a^2} + \dfrac{\xi_2^2}{b^2}+\dfrac{\xi_3^2}{c^2}=1, \quad 0<c\leq b \leq a. \nonumber
\end{equation}
We choose  $\alpha = 0.01 \, \mathrm{m}$ and  $k = 10$ and consider the cases  $a=b=c=1$ so that $B$  degenerates to a sphere and  consider the ellipsoid defined by $a = 1$, $ b = 0.7$, $c = 0.5$.  In Figures~\ref{fig:sphere:relativeError} and ~\ref{fig:ellipsoid:relativeError} we show the behaviour of $\| \mathcal{T}^{Exact} - \mathcal{T}^h\| / \|\mathcal{T}^{Exact} \|$ with uniform mesh refinement {against $1/h$} for different values of $0\le \beta \le 1$ and both primal and dual mesh discretisations for both a sphere and an ellipsoid, respectively. Here, $\| \mathcal{T}^{Exact} \| = \left ( \sum_{i,j=1}^3 ( {\mathcal T}_{ij}^{Exact)})^2 \right )^{1/2}$ is the Frobenius norm and ${\mathcal T}^{Exact}$ is the exact solution for these objects  \cite[page 83]{ammarikangbook}.
Specifically, in Figure \ref{fig:sphere:relativeError}, the primal mesh is initially chosen as surface mesh of $232$ surface triangles and then uniformly refined $6$ times until a mesh of $31\,616$ surface triangles is obtained. While in Figure \ref{fig:ellipsoid:relativeError}, the primal mesh is initially chosen as surface mesh of $48$ surface triangles and then uniformly refined $7$ times until a mesh of $65\,566$ surface triangles is obtained. 
The corresponding dual meshes are obtained from the primal meshes using the procedure described in~\cite{BetckeJCP2019}. Except for a couple of exceptions (eg  $\beta = 0.4$ shown in Figure \ref{fig:sphere:relativeError} for the dual mesh), the convergence behaviour for both primal and dual meshes is similar and practically identical for different values of $\beta$. We attribute the exceptions to some super-convergence behaviour in exceptional cases. The fixed convergence rate of the PST under uniform mesh refinement is not surprising given the results of~\ref{sec:Convergence}, which establishes that the convergence behaviour of the coefficients of the PST  is the same as the convergence of the solutions to the transmission problem under uniform mesh refinement.

\afterpage{\clearpage}

\section{Adaptive Computation of PST Coefficients using BEM++} \label{sec:AdaptiveMesh}

Apart from limited smooth geometries, such as the ellipsoid considered above, the PST does not have an analytical solution and instead numerical approximations are essential to find its coefficients. Furthermore, for realistic objects with edges and corners,  refinement of the mesh towards these features is expected to be desirable to accurately capture the solution to the transmission problem (\ref{eq:TransmissionProblem}), and, hence, the PST coefficents.  To efficiently and automatically guide this mesh refinement, we employ an existing adaptive mesh method based on \cite{BetckeJCP2019} to enhance the numerical computation ${\mathcal T}^h$ for objects with edges and corners. 
In particular, these schemes allow elements to be selected for  refinement according to the magnitude of the contributions to the error indicator with those with the largest contributions being refined. The solution on the refined mesh leading to a better approximation for an  economical increase in computational effort. 

The approach used here  closely follows the successful scheme presented for a related problem in~\cite{BetckeJCP2019} for which an existing code was available that could be suitably modified. The basic procedure could be summarised as follows:  we solve for  $\phi_{i}^h$  on a dual grid discretisation of
${\mathcal M}_0$
 and use a Zienkiewicz-Zhu type error estimator $\eta_i  \approx \| \phi_{i}^h - \phi_i \|_E$, where $\| \cdot \|_E$ denotes the energy norm. This  is obtained by interpolating solutions on the primal grid; see~\cite{BetckeJCP2019} for further details. We consider two different adaptive schemes, which are driven by flagging elements for refinement according to the elemental contributions to  $\eta = \max \eta_i$ or $\eta=\sum_{i=1}^3 \eta_i$ leading to a new mesh ${\mathcal M}_1$.  Then, by computing $\phi_{i}^h$ on ${\mathcal M}_1$, and repeating the error estimation process, a new grid ${\mathcal M}_2$ is obtained. By repeating this process, the sequence of grids
 ${\mathcal M}_\ell$, $\ell =1, \ldots$ is generated. Note that the D\"ofler criterion \cite{BetckeJCP2019} is used to mark the elements that will be refined. The refinements are controlled by $\theta$, where  $\theta=1$ coincides with uniform mesh-refinement. To maintain a balance between the number of refinements and the control of the elements to be refined, the D\"ofler parameter is chosen as $0.4 \le \theta \le 0.6$.

\section{Numerical Benchmark {Computations} of PST Coefficients using Adaptive BEM} \label{sec:numericalresults}

In this section, we present benchmark computations of the PST obtained using the aforementioned adaptive scheme.  We first define the benchmark geometries that will be considered. Then, we consider the performance in detail of the L-shape benchmark geometry before presenting benchmark {computations} for all geometries. Further computational results can be found in the open data repository~\cite{alandata}.

\subsection{Benchmark Geometries}  \label{sec:geometry}

\begin{table}[h]
\centering
\caption{Benchmark {geometries}. \label{tab:geometry:obj}} \vspace*{2pt} 
{
\setlength\arrayrulewidth{1pt}
\begin{tabular}{ c | c }
\hline \\
L-shape & Cube \\
\includegraphics[scale=0.2]{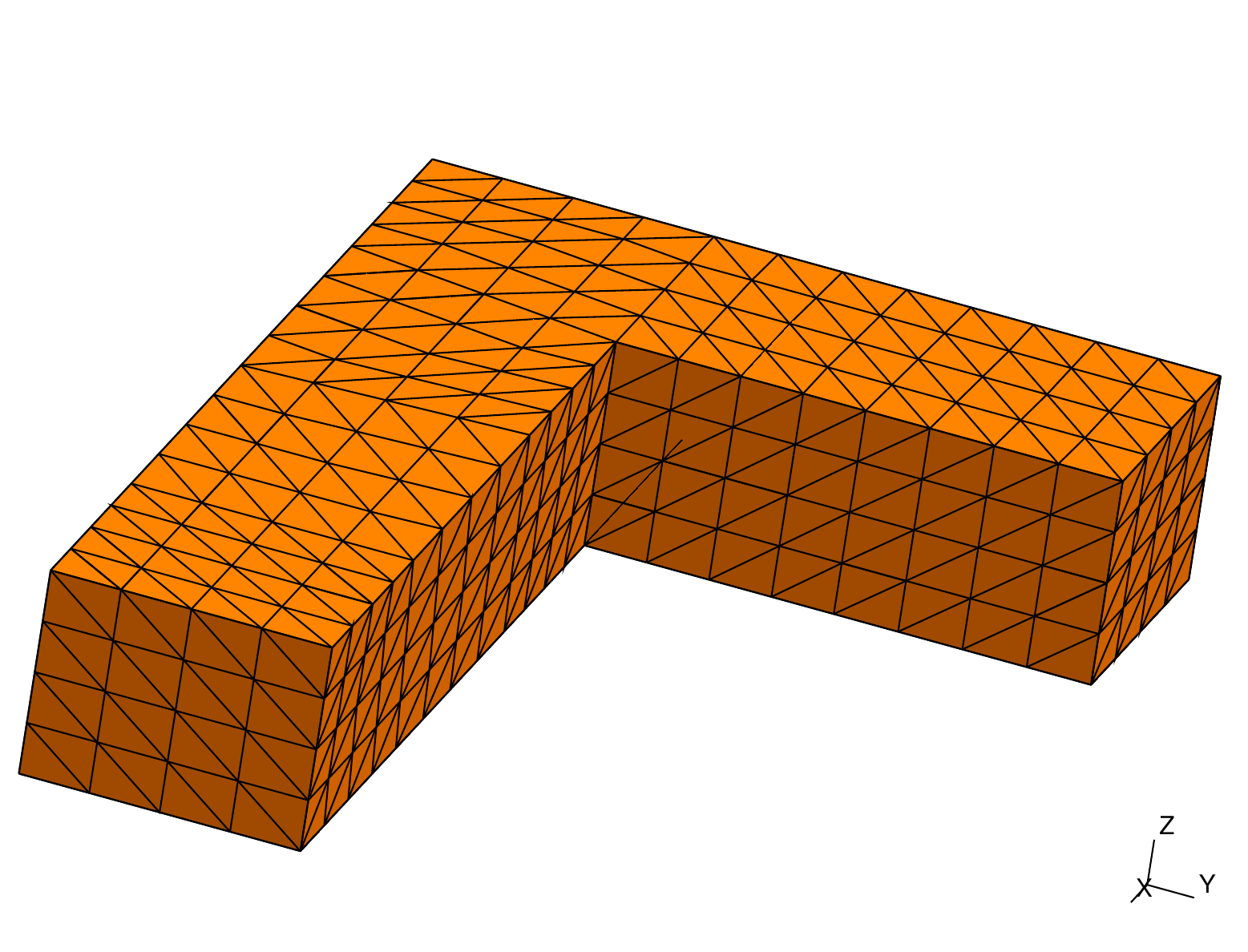}  & 
\includegraphics[scale=0.2]{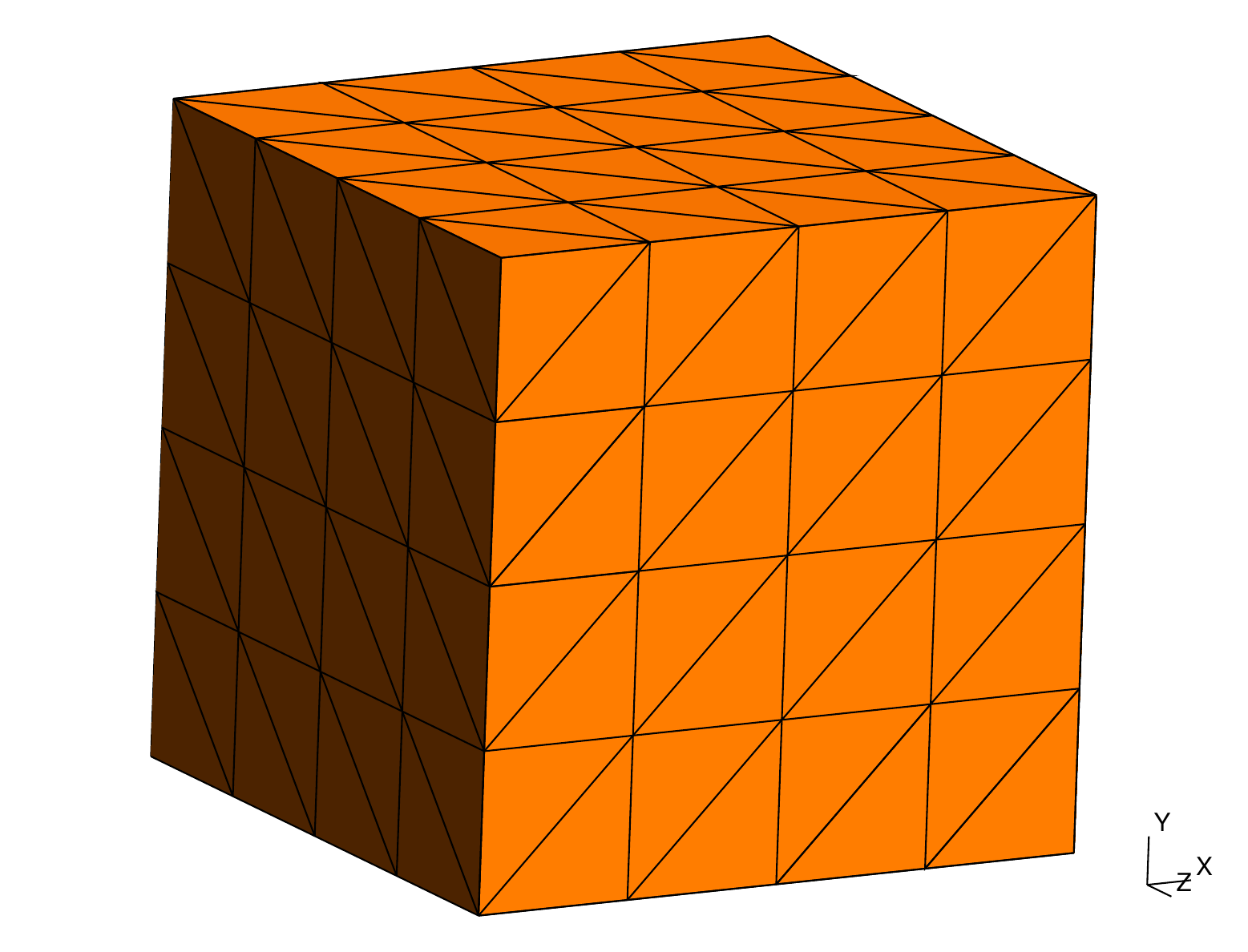}      \\
\hline\\
Tetrahedron & Key \\
\includegraphics[scale=0.2]{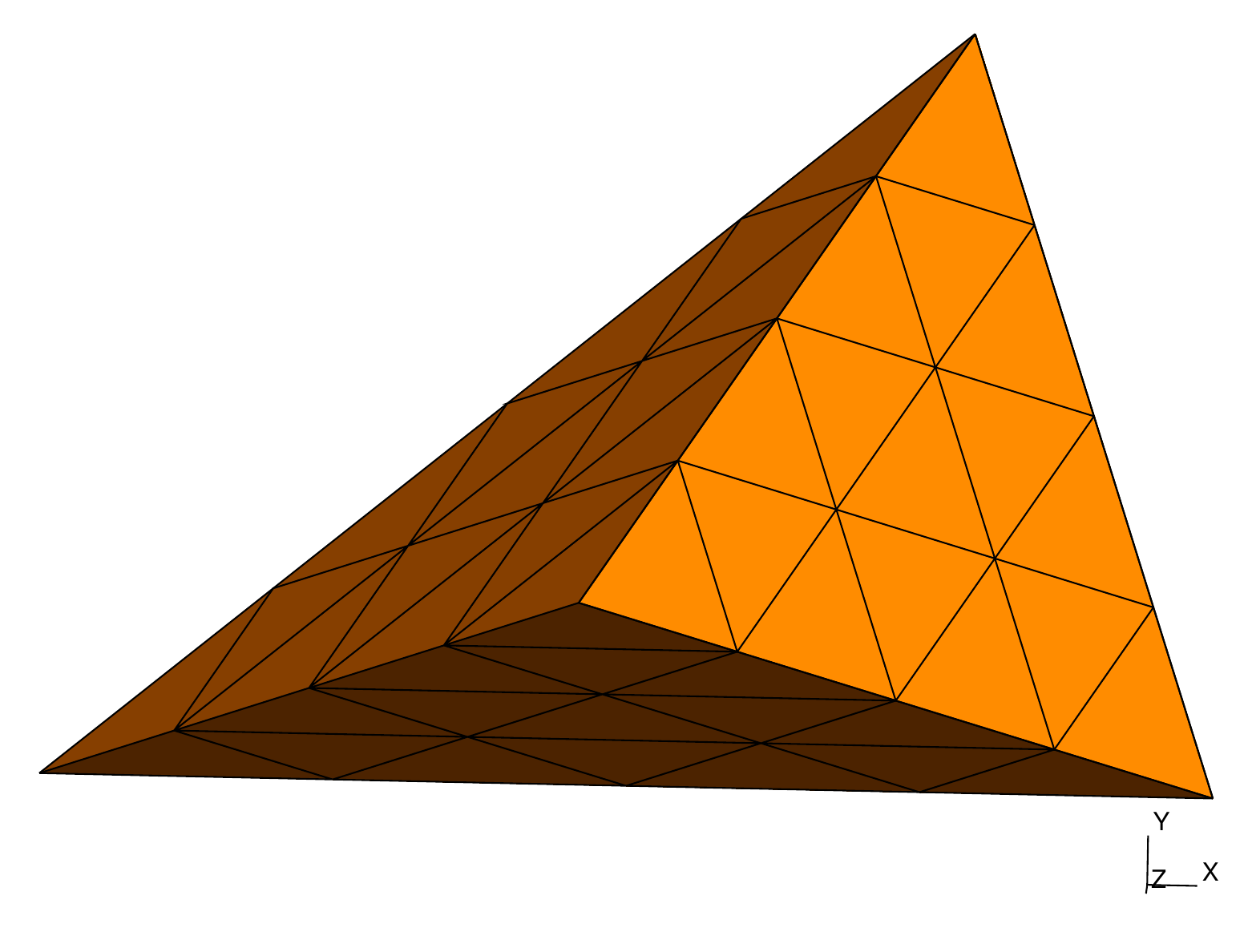} & \includegraphics[scale=0.2]{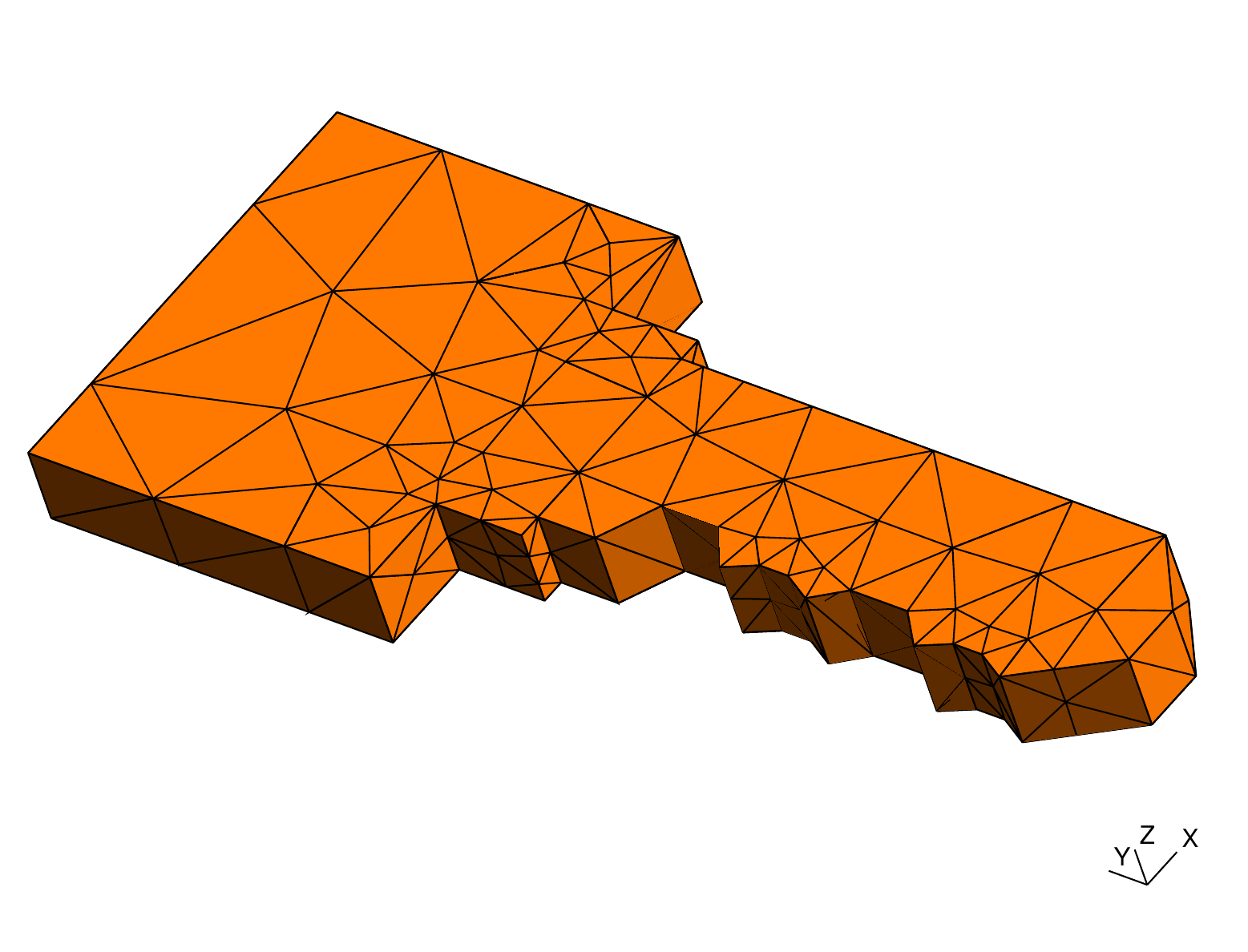}     \\
\hline
\end{tabular}
}
\end{table}

\subsubsection{L-shape}  \label{sec:geo:lshape}

For the  L-shape object, we consider $B_\alpha := \alpha B$ with $\alpha = 0.01 \, \mathrm{m}$, $B = B_1 \bigcup B_2$, where 
\begin{align*}
B_1 =& [0.0,7.8] \times [0.0,2.0] \times [0.0,1.5], \nonumber\\
B_2 = & [0.0,2.2] \times [2.0,5.6] \times [0.0,1.5]. \nonumber
\end{align*}
The reflectional symmetries of this object can be used to show that ${\mathcal T}$ {can be defined by the 4 independent coefficients ${\mathcal T}_{11}$, ${\mathcal T}_{22}$, ${\mathcal T}_{33}$  and ${\mathcal T}_{12}= {\mathcal T}_{21}$.}

\subsubsection{Cube}  \label{sec:geo:Cube}

The second object considered is where $B_\alpha := \alpha B$ is a cube, with $\alpha = 0.01 \, \mathrm{m}$ and $B$ is a unitary cube.  The rotational and reflectional symmetries of this object mean that ${\mathcal T}$ is a multiple of identity and has a single independent coefficient {${\mathcal T}_{11} = {\mathcal T}_{22} = {\mathcal T}_{33} $}.

\subsubsection{Tetrahedron}  \label{sec:geo:Tetrahedron}

The third object, we consider $B_\alpha := \alpha B$ to be a tetrahedron, with $\alpha = 0.01 \, \mathrm{m}$ and the vertices of $B$ are chosen to be at the locations 
\begin{equation}
{\vec v}_1 = \begin{pmatrix} 
0 \\
0 \\
0 
\end{pmatrix}, \ 
{\vec v}_2= \begin{pmatrix} 
7 \\
0 \\
0 
\end{pmatrix}, \  
{\vec v}_3 =  \begin{pmatrix} 
5.5 \\
4.6 \\
0.0 
\end{pmatrix}, \
{\vec v}_4=  \begin{pmatrix} 
3.3 \\
2.0 \\
5.0 
\end{pmatrix}. \nonumber
\end{equation}
This object has no symmetries and, hence, ${\mathcal T}$ has 6 independent coefficients {${\mathcal T}_{11}, {\mathcal T}_{22}, {\mathcal T}_{33}, {\mathcal T}_{12} = {\mathcal T}_{21}, {\mathcal T}_{13}={\mathcal T}_{31}$ and $ {\mathcal T}_{23}={\mathcal T}_{32}$}.

\subsubsection{Key}  \label{sec:geo:Key}

\begin{figure}[h]
\begin{center}
\includegraphics[scale=0.4]{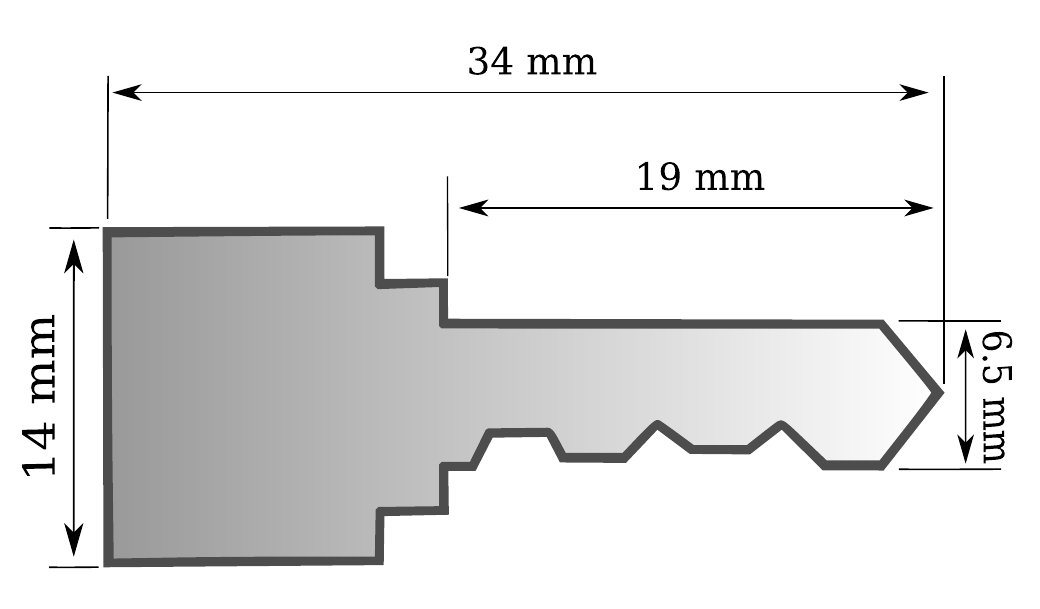}
\end{center}
  \caption{Key: Sketch of a simplified key  showing the physical object $B_\alpha$.}
        \label{fig:key:diagram}
\end{figure}

The last object, we will consider is  a  key $B_\alpha := \alpha B$, as illustrated in Figure~\ref{fig:key:diagram}, with $\alpha = 0.001 \, \mathrm{m}$, which has additional geometrical complexities. We have divided the key in three parts, bow, shoulder and blade,
\begin{equation}
B = \mathrm{bow} \bigcup \mathrm{shoulder} \bigcup \mathrm{blade}, \nonumber
\end{equation}
where 
\begin{equation*}
\mathrm{bow} = [-7.0,7.0] \times [3.0,15.0] \times [-1.25,1.25],
\end{equation*}
\begin{equation*}
\mathrm{shoulder} = [-4.0,4.0] \times {[0.0,3.0]} \times [-1.25,1.25], \qquad 
\end{equation*}
and  the maximum dimensions of the blade are
\begin{equation*}
\mathrm{blade} = [-3.25, 3.25] \times [-19.0,0.0] \times [-1.25,1.25].
\end{equation*}
To enforce that the shape is more realistic, some cuts have been made in the blade, as presented in Figure \ref{fig:key:diagram}. This object has {only a symmetry in the $\xi_3$ direction and so ${\mathcal T}$ has the 4 independent coefficients ${\mathcal T}_{11}$, ${\mathcal T}_{22}$, ${\mathcal T}_{33}$  and ${\mathcal T}_{12}={\mathcal T}_{21}$.}

\afterpage{\clearpage}
\subsection{Detailed study of the L-Shape Geometry}  \label{sec:numericalresults:Lshape}

As an illustration, we provide a detailed convergence study of the L-shape geometry. 
 In this section, we provide a convergence study as well as images of the resulting adaptive surface meshes for this object. Throughout this section, the surface of $B$ is discretised with an initial mesh ${\mathcal M}_0$ of $832$ surface triangles
and in absence of an analytical solution, we consider $\mathcal{T}^{h}_{\mathrm{Fixed}}$ obtained on a fine mesh ${\mathcal M}_{\rm{Fixed}}$ with $212\,992$ surface triangles to be indistinguishable from the exact. 

\subsubsection{Convergence Study for Varying $\beta$}  \label{sec:numericalresults:beta}

\begin{figure}[h]
\centering
    \subfigure[$\theta = 0.4$ ]{\includegraphics[scale=0.4]{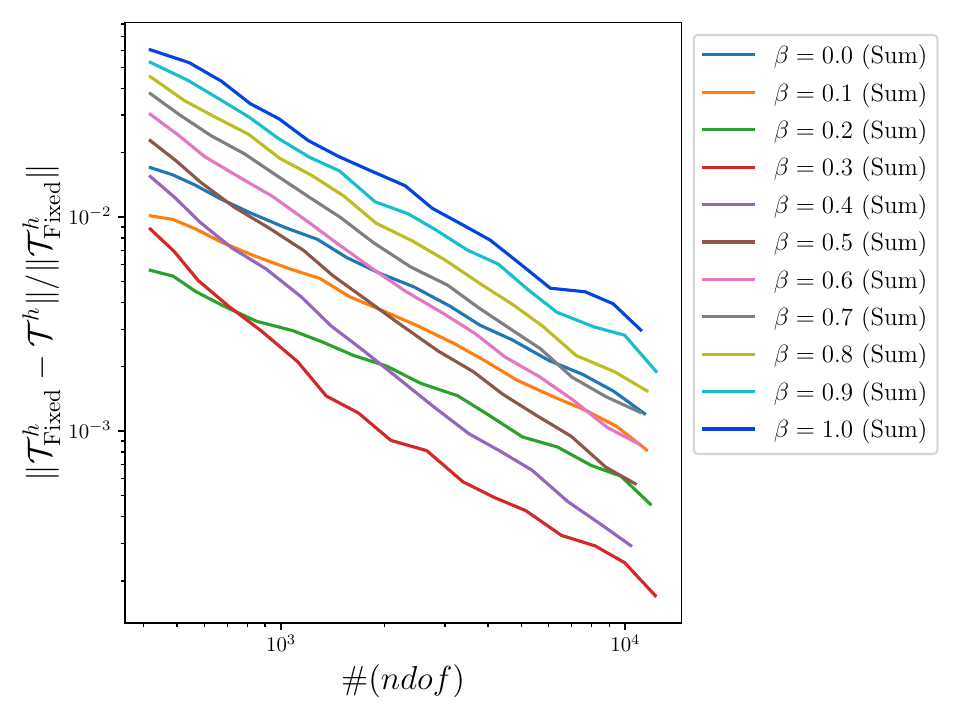} 
    \label{fig:lshape:pst:sum:4}}
    \subfigure[ $\theta = 0.5$ ]{\includegraphics[scale=0.4]{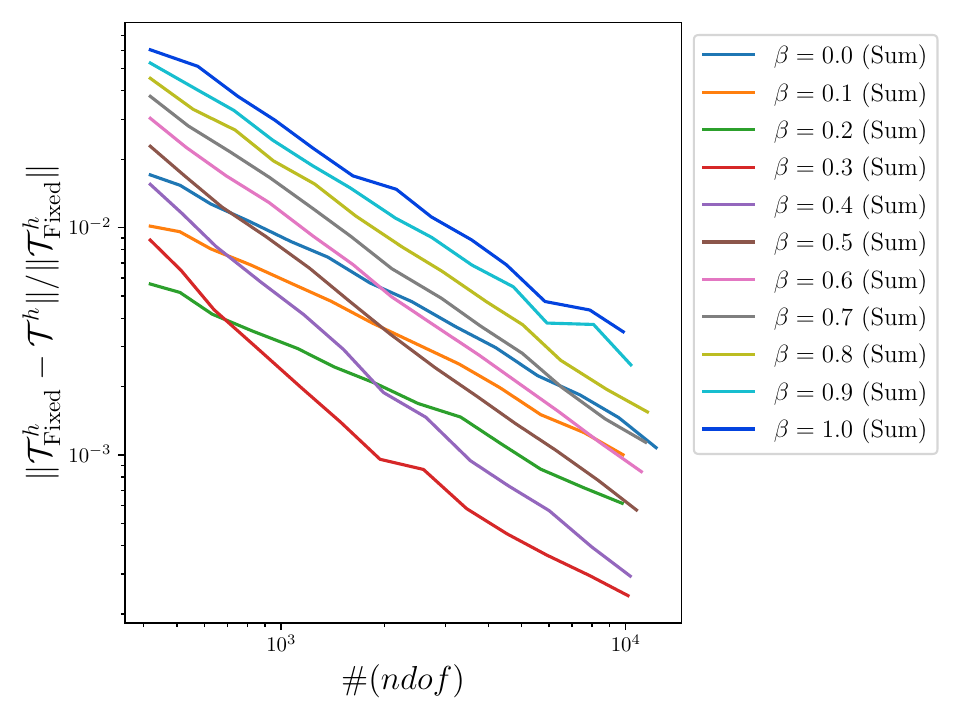} 
    \label{fig:lshape:pst:sum:5}}
     \subfigure[$\theta = 0.6$ ]{\includegraphics[scale=0.4]{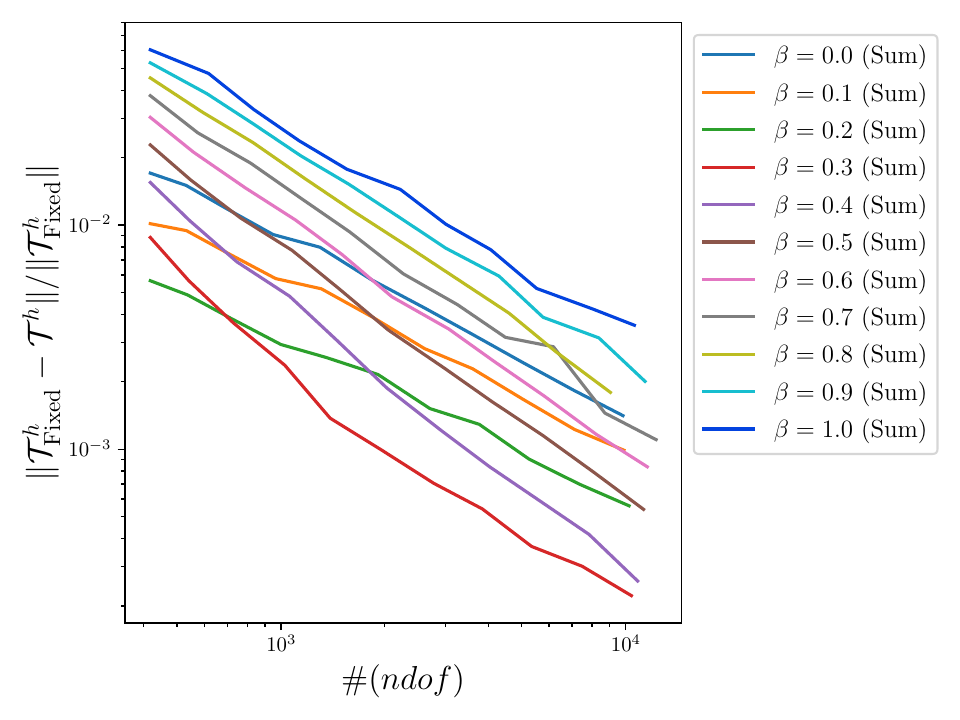} 
    \label{fig:lshape:pst:sum:6}}
  \caption{L-shape: Convergence study of ${\mathcal E} =  \left\|  \mathcal{T}^{h}_{\mathrm{Fixed}} - \mathcal{T}^{h} \right\| / \left\|  \mathcal{T}^{h}_{\mathrm{Fixed}}  \right\| $, using the $\eta= \sum_{i=1}^3 \eta_i$ estimator  for $\beta = 0.0, \ldots, 1.0$ (a) $\theta = 0.4$ (b) $\theta = 0.5$ (c) $\theta = 0.6$. }
        \label{fig:lshape:pst:sum}
\end{figure}

\begin{figure}[h]
\centering
    \subfigure[ $\theta = 0.4$  ]{\includegraphics[scale=0.4]{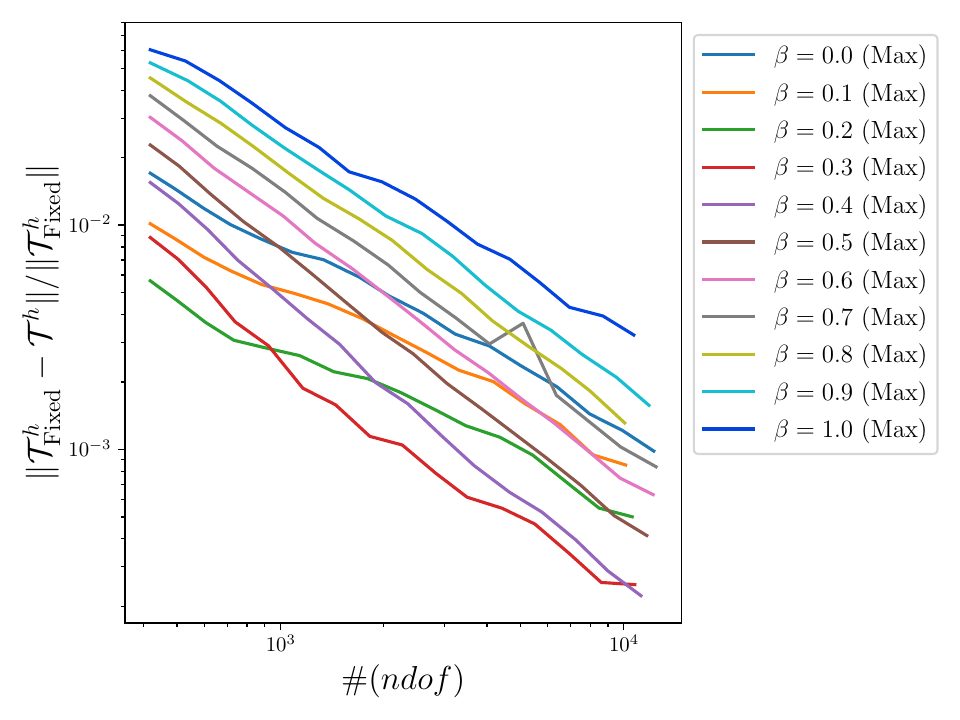} 
    \label{fig:lshape:pst:max:4}}
    \subfigure[$\theta = 0.5$ ]{\includegraphics[scale=0.4]{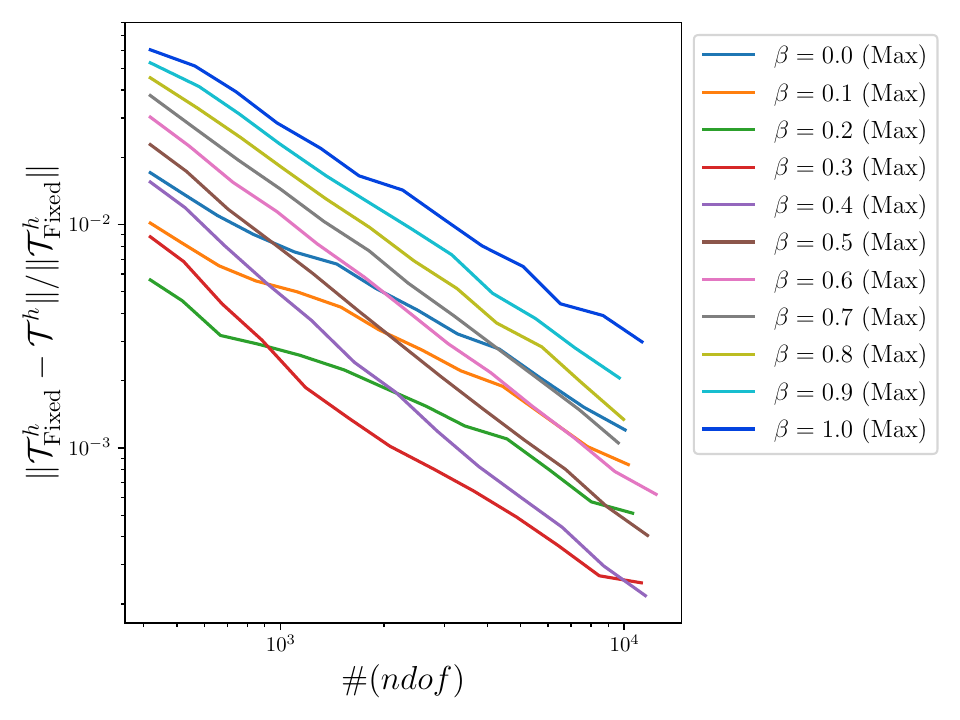} 
    \label{fig:lshape:pst:max:5}}
    \subfigure[$\theta = 0.6$ ]{\includegraphics[scale=0.4]{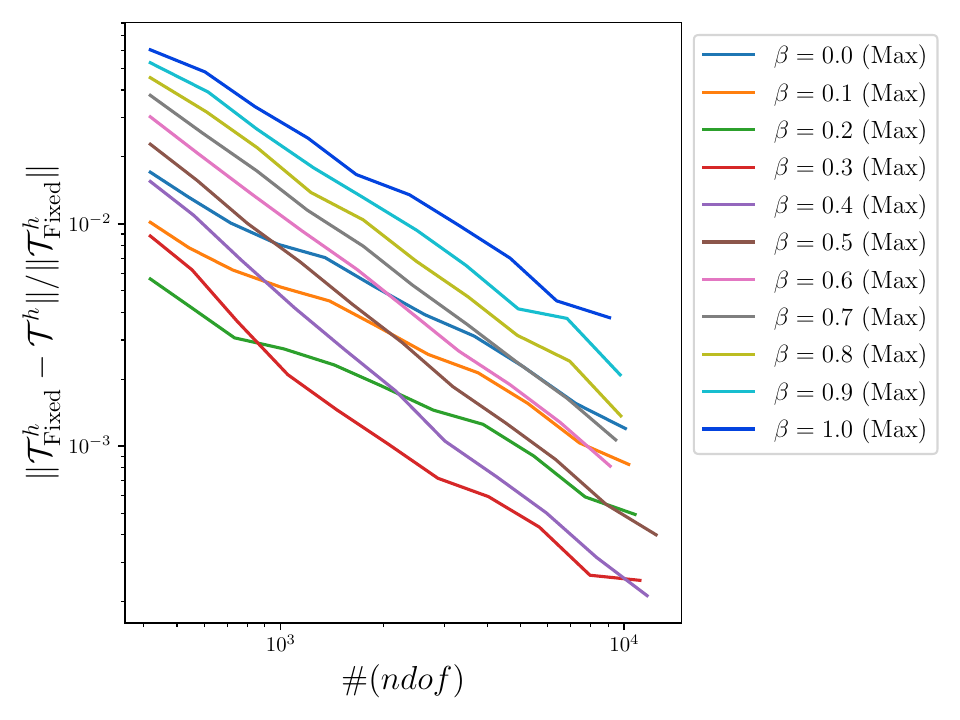} 
    \label{fig:lshape:pst:max:6}}
  \caption{L-shape: Convergence study of ${\mathcal E} =  \left\|  \mathcal{T}^{h}_{\mathrm{Fixed}} - \mathcal{T}^{h} \right\| / \left\|  \mathcal{T}^{h}_{\mathrm{Fixed}}  \right\| $, using the $\eta = \max \eta_i$  estimator  for $\beta = 0.0, \ldots, 1.0$ (a) $\theta = 0.4$ (b) $\theta = 0.5$ (c) $\theta = 0.6$. }
        \label{fig:lshape:pst:max}
\end{figure}

We  consider the accuracy of ${\mathcal T}^h$ obtained from \eqref{eq:AlternativeFormula} when  weights $\beta = 0.0, 0.1$, $\ldots, 1.0$ are considered. For a contrast $k=10$, the convergence study of the relative error of the PST given by 
\begin{equation}\label{eq:relErrorPST}
\mathcal{E} := \dfrac{\left\|  \mathcal{T}^{h}_{\mathrm{Fixed}} - \mathcal{T}^{h} \right\|}{\left\|\mathcal{T}^{h}_{\mathrm{Fixed}} \right\|},
\end{equation}
{with the number of degrees of freedom $\#(ndof)$ is} presented in Figures \ref{fig:lshape:pst:sum}--\ref{fig:lshape:pst:max}.
In Figure \ref{fig:lshape:pst:sum}, the convergence behaviour obtained when the  error estimator $\eta= \sum_{i=1}^3 \eta_i$  is used to drive the adaptive process for the cases of $\theta = 0.4$,  $\theta = 0.5$ and $\theta = 0.6$. Analogously, Figure \ref{fig:lshape:pst:max} shows the corresponding convergence behaviour obtained using  the  error estimator $\eta = \max \eta_i$.

Comparing the results shown in Figures \ref{fig:lshape:pst:sum} and \ref{fig:lshape:pst:max}, we first observe that, unlike the results presented in Figures~\ref{fig:sphere:relativeError} and~\ref{fig:ellipsoid:relativeError}, changing $\beta$ can improve the convergence rate { of ${\mathcal E}$  with respect to \#$(ndof)$. For example,  for $\theta=0.4$ and $\eta = \max \eta_i$, the convergence rate of  ${\mathcal E}$  with respect to \#$(ndof)$ can be improved from $ -1.1$ to $-2.5$ by changing $\beta=0.2$ to $\beta=0.4$ (or from $-0.75$ to $-1.4$ considering  ${\mathcal E}$  with respect to $1/h$), which can accelerate the convergence of the} adaptive scheme. In particular,
 it can be observed that for both  error estimators $\eta= \sum_{i=1}^3 \eta_i$ and $\eta = \max \eta_i$, the {the fastest convergence rate } is obtained by choosing $\beta=0.4$, which also achieves small error values {of ${\mathcal E} =  0.00021$ for the finest mesh and $\theta=0.6$}. Similar convergence curves {for ${\mathcal E}$} were found for $\beta =0.5,0.6$, but with larger errors. We have repeated this study for a range of other objects and have found that in all cases the best performance {in terms of convergence rate and accuracy for ${\mathcal E}$ were, among}  the parameters considered, for  $0.3 \le \beta \le 0.4$ and $k=10$.
{\clearpage}

\subsubsection{Convergence Study for Varying Contrasts $k$} 

We have undertaken adaptive mesh approximations for different contrasts $k= 0.1, 0.2, 1.5, 10, 15,100,150$ and summarise  the $\beta $ for which the best  results {convergence rates} were obtained using $\theta=0.6$ in Table~\ref{tab:lshape:k}.

\begin{table}[!h]
\centering
\caption{L-shape: Contrast $k$ and best $\beta$. \label{tab:lshape:k}} \vspace*{2pt} 
\begin{tabular}{ccc}
\hline 
 Contrast $k$   & Best weight $\beta$ &  Measurement option \\ 
 \hline
$0.1$     &  $0.0$  &    max   \\
$0.2$     &  $0.4$  &    max   \\
$1.5$     &  $0.6$  &    max   \\
$10.0$   &  $0.4$  &    max   \\
$15.0$   &  $0.3$  &    sum   \\
$100.0$ &  $0.1$  &    sum   \\
$150.0$ &  $0.1$  &    sum   \\
\hline
\end{tabular}
\end{table}

\begin{figure}[h]
\centering
    \subfigure[Comparison of the relative error for $k=0.2$]{\includegraphics[scale=0.4]{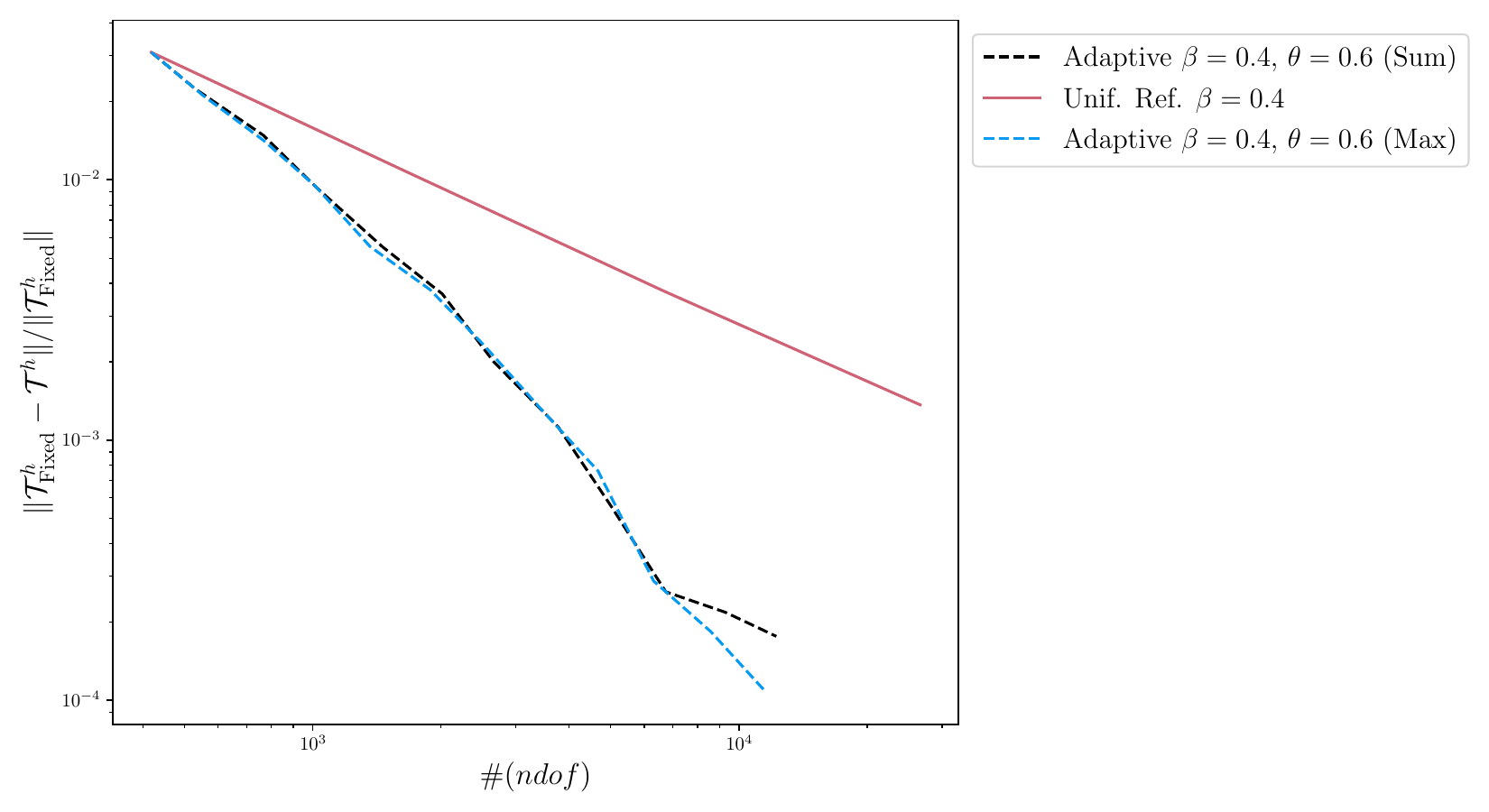} 
    \label{fig:lshape:relErrorPST:k02}}
     \subfigure[ZZ-error estimator for $k=0.2$]{\includegraphics[scale=0.4]{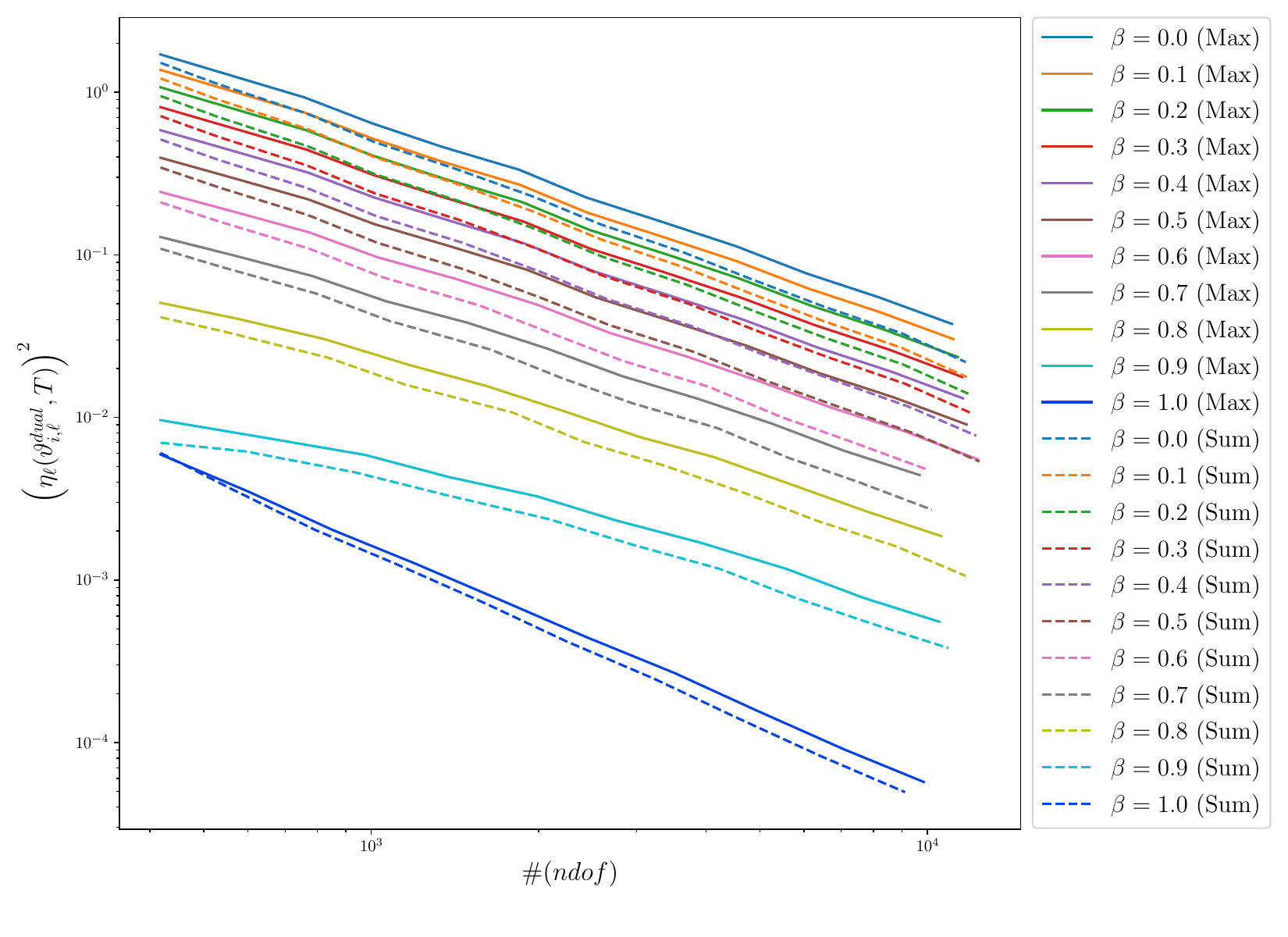} 
    \label{fig:lshape:eta:k02}}
  \caption{L-shape: (a) Comparative study between uniform refinement and adaptive mesh solutions for ${\mathcal E} =  \left\|  \mathcal{T}^{h}_{\mathrm{Fixed}} - \mathcal{T}^{h} \right\| / \left\|  \mathcal{T}^{h}_{\mathrm{Fixed}}  \right\| $ for $\beta = 0.4$, $\theta = 0.6$ and $k=0.2$ (b) Convergence study of the ZZ-error estimator, using  $\eta= \sum_{i=1}^3 \eta_i$ and $\eta = \max \eta_i$  for $\beta = 0.0, \ldots, 1.0$, $\theta = 0.6$ and $k=0.2$. }
        \label{fig:lshape:k02}
\end{figure}

If adaptive mesh approximations for different contrasts $k= 0.1, 0.2, 1.5, 10$, $15,100,150$ are compared with uniform mesh refinements, we find that the adaptive algorithm offers a superior performance {in terms of convergence rate and accuracy for ${\mathcal E} $ compared to} uniform refinement for $ k= 0.1, 0.2, 1.5, 10, 15$ with the best performance
for $k=0.1,0.2$. For $k=100, 150$, the adaptive approach does not show an improvement over uniform refinement. We illustrate one of the best performing cases {in terms of improvement in convergence rate and accuracy} corresponding to $k=0.2$ in Figure~\ref{fig:lshape:relErrorPST:k02} and show how error estimators  $\eta= \sum_{i=1}^3 \eta_i$ and $\eta = \max \eta_i$ converge for this case and different values of $\beta$ in Figure~\ref{fig:lshape:eta:k02}. In these experiments, a maximum of $25\,000$ surface triangles was generated by the adaptive algorithm while the uniform refinement led to $65\,536$ surface triangles.

\subsubsection{Illustrations of Adaptive Meshes} 

\begin{figure}[h]
    \centering
    \subfigure[Max - $\ell = 0$ ]{\includegraphics[scale=0.2]{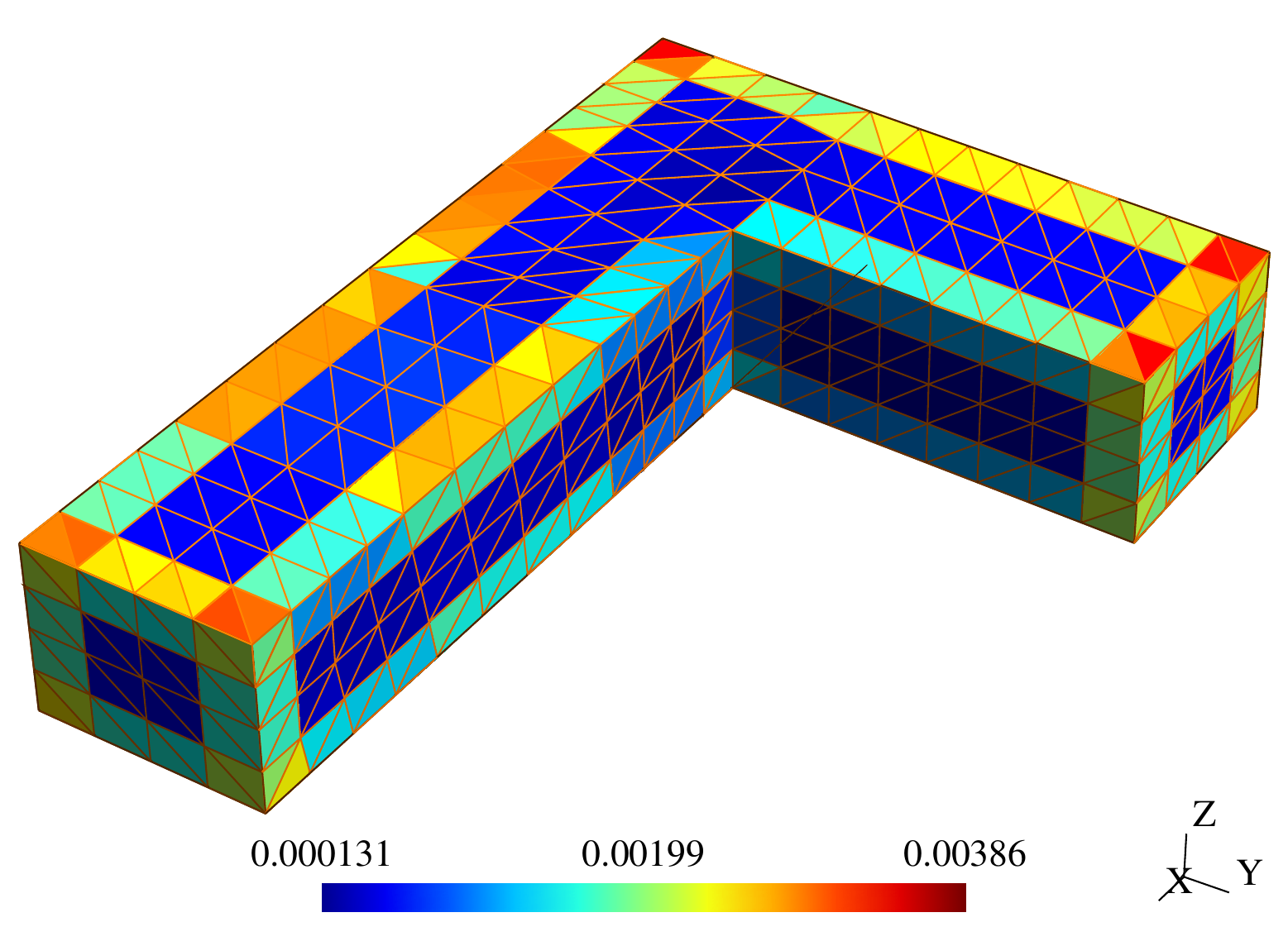} 
    \label{fig:lshape:max:0}}
    \subfigure[Max - $\ell = 2$ ]{\includegraphics[scale=0.2]{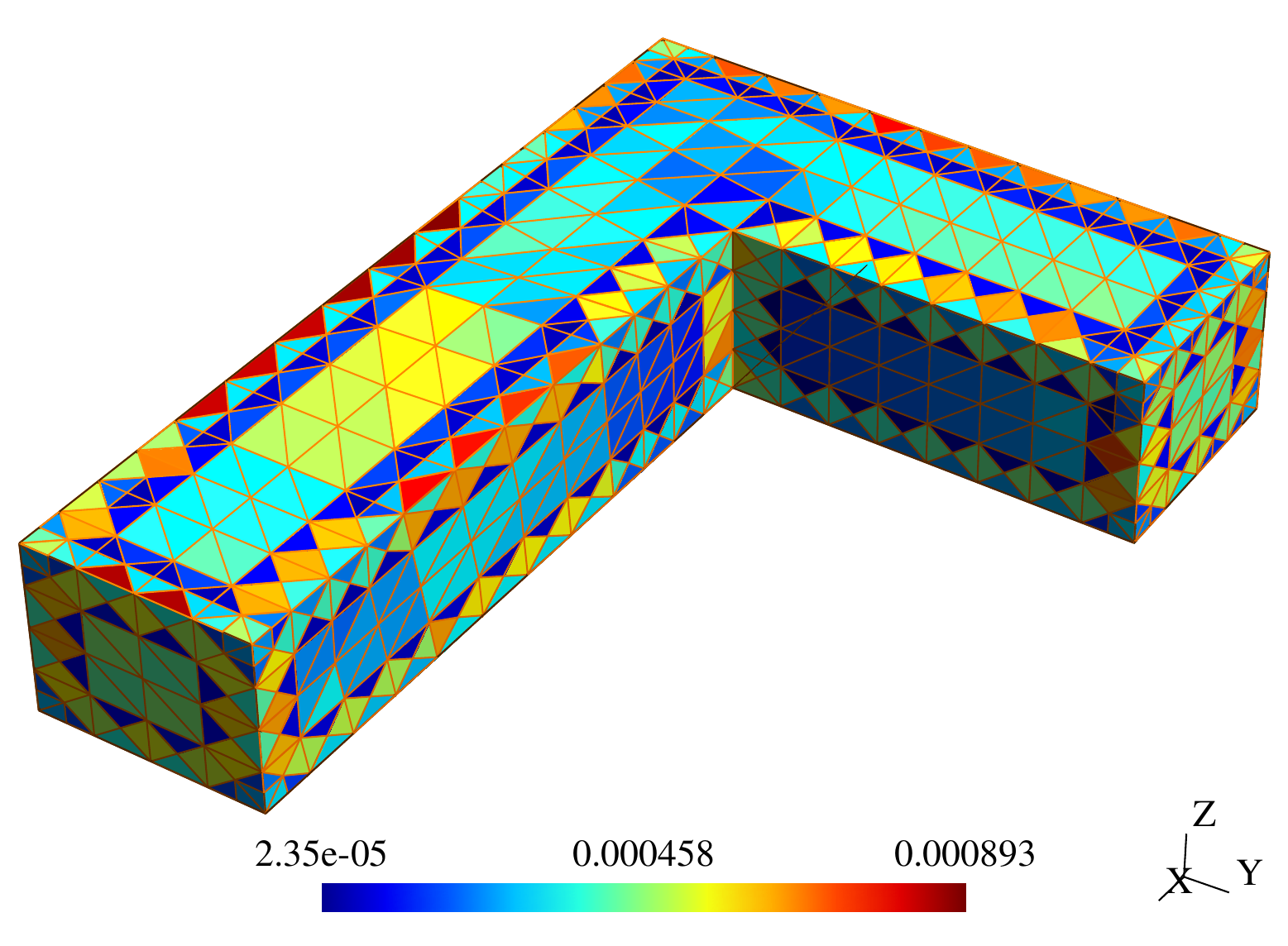} 
    \label{fig:lshape:max:2}}
    \subfigure[Max - $\ell = 4$ ]{\includegraphics[scale=0.2]{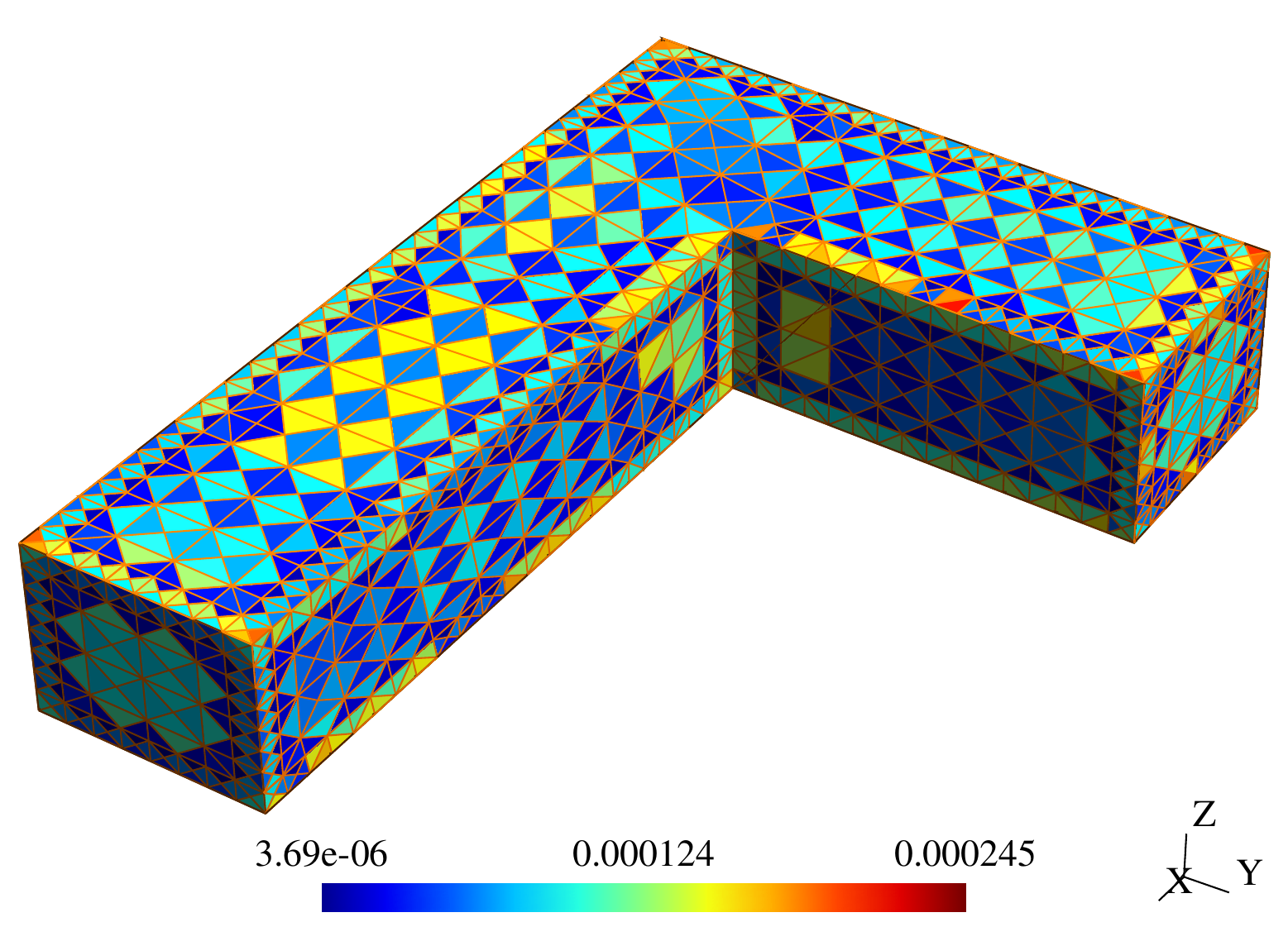} 
    \label{fig:lshape:max:4}}
    \subfigure[Max - $\ell = 6$ ]{\includegraphics[scale=0.2]{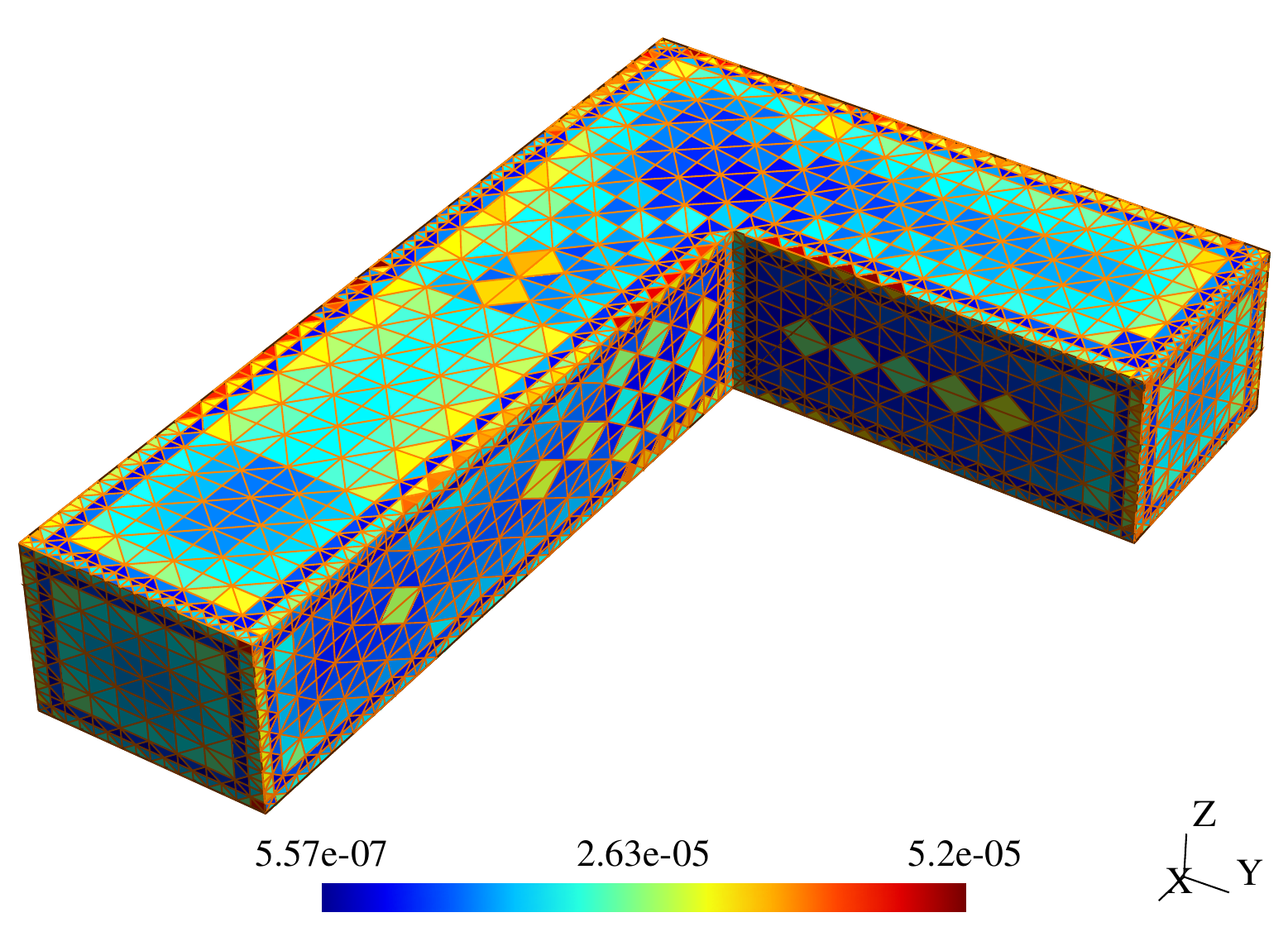} 
    \label{fig:lshape:max:6}}
    \caption{L-shape: Adaptive mesh generated using the ZZ-type error estimator using {\it max} direction measurement $\eta_\ell( \vartheta_i^{dual},T)$ for $\mathcal{M}^{dual}_\ell$ (a) $\ell = 0$ (b) $\ell = 2$ (c) $\ell = 4$ (d) $\ell = 6$.}
    \label{fig:lshape:adaptiveMesh:max}
\end{figure}

For the case of $k=10$, the initial grid ${\mathcal M}_0$ and exemplar adaptive grids  ${\mathcal M}_\ell$, $\ell =2, 4, 6$ obtained for $\beta=0.6$, $\theta=0.6$ using   $\eta = \max \eta_i$ are shown in Figures \ref{fig:lshape:max:0}-\ref{fig:lshape:max:6},  respectively, which have $832$ to $6\,022$ surface triangles.
Note that the meshes generated by the algorithm are more refined along the edges and corners, this trend is expected, since the solution $\phi_{i,h}$ is less regular at these locations and elements in the vicinity of edges and corners contribute most to the error estimate $\eta_i$. The distribution of the refined elements are similar for other objects.

\clearpage
\subsubsection{Comparison of Computational Time} 
{As an example of computational timings, the L-shape object with contrast $k=0.2$, weight $\beta = 0.4$, and $\theta = 0.6$ {and a target accuracy of $\mathcal{E} \approx 1.3 \times 10^{-3}$ is considered. The time taken by the adaptive mesh and the uniform mesh approaches  is compared in} Table \ref{tab:computationalTime}, where a time saving {of around 9 minutes} can be observed. For this example, we have used the following machine {\it Intel(R) Xeon(R) CPU E5-2687W v2 at 3.40GHz, 8 cores} and {\it 128 GB of RAM}.}

\begin{table}[htb]
\centering
\caption{L-shape: Computational time study. \label{tab:computationalTime}} \vspace*{2pt} 
\begin{tabular}{cccc}
\hline
Refinement approach   & $\mathcal{E}$         &  $\#ndof$    &  Computational Time\\ 
\hline
Uniform Mesh     &  $1.376036530 \times 10^{-3}$  &    $26\,626$ &    0:14:22   \\
Adaptive Mesh     & $1.330921701 \times 10^{-3}$  &    $3\,412$  &    0:05:12   \\
\hline
\end{tabular}
\end{table}

\subsection{Reference Polarization Tensors and Benchmark {Computations}}  \label{sec:reference:PST}

The benchmark {computations} obtained using the adaptive mesh algorithm for each of the benchmark geometries in Section~\ref{sec:geometry} are now presented for the case of $k=10$. We state ${\mathcal E} $ for each of the tensors resulting from the final adaptive mesh. In addition, by introducing the splitting
\begin{equation}
{\mathcal T} = {\mathcal O} + {\mathcal D}, \nonumber
\end{equation}
where ${\mathcal D}$ is the diagonal part of ${\mathcal T}$ and ${\mathcal O}$ the off-diagonal part,
we provide a measure of the error in the off-diagonal part of the tensor as
\begin{equation}
\mathcal{E}^{\mathrm{off}} := \dfrac{\left\|  \mathcal{O}^{h}_{\mathrm{Fixed}} - \mathcal{O}^{h} \right\|}{\left\|\mathcal{T}^{h}_{\mathrm{Fixed}} \right\|}. \nonumber
\end{equation}
for $\mathcal{T}^h$ obtained on the final adaptive mesh. The results are summarised in Table \ref{tab:best:error} for each benchmark {geometry}, using $\eta = \max \eta_i$ and $\theta = 0.6$.

\begin{table}[h]
\centering
\caption{Relative error and off-diagonal error for the best convergence result. \label{tab:best:error}} 
\begin{tabular}{c|c|c|c|c|c}
\hline 
Object    &  Best $\beta$    &  Maximum&  Number of elements      &  $\mathcal{E}$     &   $\mathcal{E}^{\mathrm{off}}$    \\
 {}  & {} &  $\#ndof$  & on $\mathcal{M}$  & {in \%} & {in \%} \\
\hline 
L-shape            &  $0.4$   &  $11\,681$     &  $23\,358$  & $0.0213$     &    0.0019 \\
Cube                &   $0.4$   &  $11\,522$    &  $23\,040$  & $0.0013$     &    0.0001 \\
Tetrahedron     &   $0.5$   &  $12\,078$    &  $24\,152$  & $0.0170$     &    0.0140 \\
Key                  &   $0.3$   &  $17\,179$    &  $34\,354$  & $0.0132$     &    0.0051 \\
\hline
\end{tabular}
\end{table}

In the following sections, we present for each benchmark {geometry}  the solution obtained with the {smallest ${\mathcal E}$
  from the adaptive mesh method} with the PST coefficients arranged as a $3 \times 3$ matrix in each case.  Note that the solutions are restricted to a mesh with a maximum of $25\,000$ surface triangles, the only exception being the key object that the mesh can contain up to $40\,000$ surface triangles. We believe that at least the first four significant figures in our benchmark tensor characterisations are accurate.

\subsubsection{L-shape: Best Approximation}  

The best approximation for the L-shape object is obtained  for $\beta=0.4$, $\theta=0.6$ and $\eta = \max \eta_i$ and has a final adapted mesh with $23\, 358$ surface triangles. This approximate tensor $\mathcal{T}^{h}$ is shown in \eqref{eq:lshape:best:tensor}
\begin{align} \label{eq:lshape:best:tensor}
\frac{\mathcal{T}^{h}}{10^{-4}} =  & \left ( \begin{array}{ccc} 
\phantom{-}1.54683606\phantom{\times 10^{-6}} &        
-0.15091553\phantom{\times 10^{-6}}                             
-6.99198691\times 10^{-6}     \\                                         
-0.15091553\phantom{\times 10^{-6}}&                          
\phantom{-}1.16202896\phantom{\times 10^{-6}}         
\phantom{-}3.33569676\times 10^{-6}\\                             
-6.99198691\times 10^{-6}&                                              
\phantom{-}3.33569676\times 10^{-6}                             
\phantom{-}0.57033981\phantom{\times 10^{-6}}            
\end{array} \right ).
\end{align}

\subsubsection{Cube: Best Approximation}  

For the cube object the best approximation is obtained by choosing $\beta=0.4$, $\theta=0.6$ and $\eta = \max \eta_i$ on a mesh with $23\, 040$ surface triangles. This approximate tensor $\mathcal{T}^{h}$ is shown in \eqref{eq:cube:best:tensor}
\begin{align} \label{eq:cube:best:tensor}
\frac{\mathcal{T}^{h}}{10^{-6}} =   & \left ( \begin{array}{ccc}  
\phantom{-}2.51110996\phantom{\times 10^{-6}} &        
\phantom{-}1.92945670\times 10^{-7} &                          
-1.59571932\times 10^{-6} \\                                             
\phantom{-}1.92945670\times 10^{-7}&                          
\phantom{-}2.51111340\phantom{\times 10^{-6}}  &       
-4.11272795\times 10^{-6}\\                             
-1.59571932\times 10^{-6}&                                              
-4.11272795\times 10^{-6}         &                    
\phantom{-}2.51110887\phantom{\times 10^{-6}}            
\end{array} \right ).
\end{align}

\subsubsection{Tetrahedron: Best Approximation}  

For the tetrahedron the best approximation is obtained by choosing $\beta=0.5$, $\theta=0.6$ and   $\eta = \max \eta_i$  on a mesh with $24\, 152$ surface triangles. This approximate tensor $\mathcal{T}^{h}$ is shown in \eqref{eq:tetra:best:tensor}

\begin{equation} \label{eq:tetra:best:tensor}
\mathcal{T}^{h} = 10^{-5} \begin{pmatrix} 
\phantom{-}9.30682676 &        
1.12847255&                          
-0.76359289 \\                                             
\phantom{-}1.12847255&                          
6.83952305  &                         
\phantom{-}0.43032199\\                             
-0.76359289&                                              
0.43032199 &                            
\phantom{-}7.80618516           
\end{pmatrix}.
\end{equation}
	
\subsubsection{Key: Best Approximation}  

For the key the best approximation is obtained by choosing $\beta=0.3$, $\theta=0.6$  and $\eta = \max \eta_i$ on a mesh with $34\, 354$ surface triangles. This approximate tensor $\mathcal{T}^{h}$ is shown in \eqref{eq:key:best:tensor}.
\begin{align}
\label{eq:key:best:tensor}
\frac{\mathcal{T}^{h}}{10^{-6}} =  & \left ( \begin{array}{ccc} 
\phantom{-}2.66099087\phantom{\times10^{-7}} &  
-2.40425167\times 10^{-2} &                      
\phantom{-}1.94050498\times 10^{-5} \\                                       
-2.40425167\times 10^{-2} &                     
\phantom{-}4.49014335\phantom{\times10^{-5}} & 
-3.78250864\times 10^{-6} \\                    
\phantom{-}1.94050498\times 10^{-5} &                                      
-3.78250864\times 10^{-6}         &           
\phantom{-}0.96008445\phantom{\times10^{-7}}    
\end{array} \right ).
\end{align}

\section{Conclusions} \label{sec:Conclusions}

{In this work, we have proposed a series of benchmark {computations} for the P\'oyla-Szeg\"o tensor (PST) for different objects.  We expect these benchmark {computations} to be of interest to electrical impedance tomography (EIT) practitioners, and, in particular, for those developing tools that will form part of a machine learning (ML) classification algorithm for detecting and classifying small conducting inclusions,  but also to other computational partial differential equation solver developers.}

{PST object characterisations are attractive as they can easily be determined from voltage measurements in EIT, once the inclusion is known, by solving an over-determined linear system for their coefficients following from (\ref{eq:asymptoticExpansion}). This has advantages over a using voxelated grid and solving 
 an ill-posed inverse problem for conductivity values in each voxel, which requires careful regularisation. However, in order for classification approaches based on PST object descriptions to work effectively, it is important that an accurate tool is used to compute the  PST characterisations and our benchmark computations allow software developers to check that this is indeed the case. Nonetheless, a limitation of the PST object characterisation is that it only characterises an object up to the best fitting ellipsoid and does not allow the conductivity and the shape of the inclusion to be uniquely determined, which ultimately limits classification approaches built on PST descriptions. This can be improved by taking advantages of the spectral behaviour of PST coefficients by taking voltage measurements as a function of frequency, or, alternatively, by considering generalised polarizability tensor object characterisations, which include more information about an object's shape and its materials.}

{We have presented a series of alternative boundary integral formulations, which, although identical for exact arithmetic, differ in practical computations for the computation of PST coefficients. We  have described how these formulations can implemented in the boundary element method python package (Bempp). We have discussed {how the application of an existing adaptive algorithm automates the refinement of the boundary element grid and} reduces the computational effort required to compute the PST coefficients to a desired level of accuracy. We have included a series of numerical examples to demonstrate these benefits for  objects with sharp corners and edges such as an L-shape domain, a cube, a tetrahedron and a key and explored how the conductivity constant effects the performance of the adaptive algorithm. }

\appendix

\section{Convergence of the P\'olya-Szeg\"o Tensor} \label{sec:Convergence}

In this section, a convergence result for the numerical approximation of ${\mathcal T}^h$ to  ${\mathcal T}$ under uniform mesh refinement is presented. 
  
  In line with the solution to related boundary integral approximations of scalar Laplace transmission problems (e.g.~\cite[pg. 199, Proposition 4.1.31]{Sauterbook2010}), we conjecture that the approximate solution $\phi_i^h$ to $\phi_i$ satisfies the following a-priori type error estimate under uniform mesh refinement
    \begin{equation}\label{eq:lemma:Sauter:3}
        \left\| \phi_i- \phi_i^h \right\|_{L^2(\Gamma)} \leq C h^s \| \phi_i \|_{H^s(\Gamma)}, 0 \leq s \leq 1.
    \end{equation}

It follows that we can establish that the convergence of the coefficients of the PST also have the same rate
    \begin{equation*}
        \left\| {\mathcal T} - {\mathcal T}^h \right\| \leq C h^s, \quad 0 \leq s \leq 1.
    \end{equation*}

To show this, consider the coefficients of the exact PST ${\mathcal T}$, using the LP form, which are given by
    \begin{equation*}
        {\mathcal T}_{ij} =  \int_{\Gamma} \xi_j \phi_i\,\dif \bm{ \xi}.
    \end{equation*}
    Let ${\mathcal T}^h_{ij}$ be a tensor computed numerically written as 
    \begin{equation*}
        {\mathcal T}^h_{ij} =  \int_{\Gamma} \xi_j \phi_{i}^h \,\dif \bm{ \xi},
    \end{equation*}
    where $\phi_{i}^h$ is an approximate solution computed using \eqref{eq:BEM:LP} and the effects of approximate numerical integration are ignored. 
    Thus, 
    \begin{align}
        \left\| {\mathcal T} - {\mathcal T}^h \right\|
        =& \left( \sum^3_{i,j=1} \left| {\mathcal T}_{ij} - {\mathcal T}^h_{ij} \right|^2 \right)^{1/2} \nonumber \\ 
        =& \left( \sum^3_{i,j=1} \left|  \int_{\Gamma} \xi_j \left( \phi_i - \phi_{i}^h \right) \,\dif \bm{ \xi} \right|^2 \right)^{1/2} \nonumber\\ 
        \leq& \left( \sum^3_{i=1} \left|  \int_{\Gamma} \left( \phi_i - \phi_{i}^h \right) \,\dif \bm{ \xi} \right|^2 \right)^{1/2} \left( \sum^3_{j=1} \left|  \int_{\Gamma} \xi_j  \,\dif \bm{ \xi} \right|^2 \right)^{1/2},  \nonumber 
        \end{align}
where the Cauchy-Schwarz inequality has been applied. This can be written as 
        \begin{align*}
       \le& \left ( \sum_{i=1}^3 \left \| \phi_i - \phi_{i}^h \right \|_{L^2(\Gamma)}^2 \right)^{1/2}  \left ( \sum_{j=1}^3 \| \xi_j \|_{L^2(\Gamma)}^2 \right)^{1/2}.
    \end{align*}
 Note that, $ \|   \xi_j   \|_{L^2(\Gamma)}   \leq C$
    with $C$ independent of $h$ and
then we have
    \begin{align}
        \left\| {\mathcal T} - {\mathcal T}^h \right\|
        \leq& C \left( \sum^3_{i=1} \left\|  \phi_i - \phi_{i}^h \right\|_{L^2(\Gamma)}^2 \right)^{1/2}.
        \label{eq:proof:estimate:1}
    \end{align}
    Substituting the error estimate in  \eqref{eq:lemma:Sauter:3} into \eqref{eq:proof:estimate:1}, we have
    \begin{align*}
        \left\| {\mathcal T} - {\mathcal T}^h \right\|
        \leq&  C h^s \max_{i=1,2,3} \| \phi_i \|_{H^s(\Gamma)},
    \end{align*}
    with $C$ independent of $h$. Note that,
    \begin{align*}
        \| \phi_i \|_{H^s(\Gamma)}
        \leq&  C,
    \end{align*}
    where the constant $C$  depends of the shape of $B$, but not on $h$. Hence leading to the quoted result.
\section*{Acknowledgement}
{A.A.S. Amad and P.D. Ledger
 gratefully acknowledges the financial support received from EPSRC in the
form of grant EP/R002134/2 and T. Betcke gratefully acknowledges the financial support received from EPSRC in the form of grant  EP/R002274/1. D. Praetorius acknowledges the support of the Austria Science Fund FWF
through grant P27005 and through the special research program {\em Taming complexity in partial differential systems} (grant F65).}

\bibliographystyle{elsarticle-num}
\bibliography{bibdatabase}

\end{document}